\documentclass{amsart}

\usepackage
[left=1.2in, right=1.2in, top=1in, bottom=1in, marginparwidth=1.0in]
{geometry}

\usepackage{amsmath,amsfonts,amsthm,amssymb}
\usepackage{array}
\usepackage[british]{babel}
\usepackage{bbm}
\usepackage{booktabs}
\usepackage{comment}
\usepackage{enumitem}
\usepackage[mathscr]{euscript}
\usepackage{float}
\usepackage[T1]{fontenc}
\usepackage{graphicx}
\usepackage{hyperref}
\usepackage[utf8]{inputenc}
\usepackage{longtable}
\usepackage{lscape}
\usepackage{mathrsfs}
\usepackage{mathtools}
\usepackage[framemethod=tikz]{mdframed}
\usepackage{multicol}
\usepackage{multirow}
\usepackage{pbox}
\usepackage{pict2e}
\usepackage{relsize}
\usepackage{soul}
\usepackage{stmaryrd}
\usepackage{tabu}
\usepackage{tikz}
\usepackage{tikz-cd} 
\usepackage{todonotes}
\usepackage[usegeometry,paper=8.5in:11in]{typearea}
\usepackage{verbatim}
\usepackage[normalem]{ulem}
\usepackage[all]{xy}

\KOMAoptions{twoside=false}

\usepackage[noabbrev,capitalise]{cleveref}
\hypersetup{%
  bookmarksnumbered=true,%
  colorlinks=true,%
  linkcolor=blue,%
  citecolor=blue,%
  filecolor=blue,%
  menucolor=blue,%
  urlcolor=blue,%
  pdfnewwindow=true,%
  pdfstartview=FitBH}   

\numberwithin{equation}{section}
\newtheorem{thm}[equation]{Theorem}
\newtheorem{cor}[equation]{Corollary}
\newtheorem{prop}[equation]{Proposition}
\newtheorem{lem}[equation]{Lemma}
\newtheorem{Que}[equation]{Question}
\theoremstyle{definition}

\newtheorem{ex}[equation]{Example}
\newtheorem{problem}[equation]{Problem}
\newtheorem{rem}[equation]{Remark}
\newtheorem{conj}[equation]{Conjecture}

\newcommand{\mmf}{\mathit{mmf}}
\newcommand{\tmf}{\mathit{tmf}}

\DeclareMathOperator{\im}{im}
\newcommand{\C}{\mathbb{C}}
\newcommand{\Z}{\mathbb{Z}}
\newcommand{\F}{\mathbb{F}}
\newcommand{\N}{\mathbb{N}}
\newcommand{\kappabar}{\overline{\kappa}}
\newcommand{\cA}{\mathcal{A}}
\newcommand{\cB}{\mathcal{B}}
\newcommand{\cP}{\mathcal{P}}
\newcommand{\cl}{\mathrm{cl}}
\DeclareMathOperator{\Ext}{Ext}
\DeclareMathOperator{\coker}{coker}

\definecolor{gray}{gray}{0.4}
\definecolor{green}{rgb}{0, 0.65, 0}
\definecolor{purple}{rgb}{0.67, 0, 1}

\begin{document}

\title[Exotic periodic phenomena]{Exotic periodic phenomena in the cohomology of the moduli stack of $1$-dimensional formal group laws}

\author[Isaksen]{Dan Isaksen}
\address{Wayne State University}
\email{isaksen.dan@gmail.com}

\author[Kong]{Hana Jia Kong}
\address{School of Mathematical Sciences, Zhejiang University, Hangzhou, China}\email{hjkong@zju.edu.cn}

\author[Li]{Guchuan Li}
\address{School of Mathematical Sciences, Peking University, Beijing, China}\email{liguchuan@math.pku.edu.cn}

\author[Ruan]{Yangyang Ruan}
\address{Beijing Key Laboratory of Topological Statistics and Applications for Complex Systems, Beijing Institute of Mathematical Sciences and Applications, Beijing 101408, China}\email{ruanyy@amss.ac.cn}

\author[Zhu]{Heyi Zhu}
\address{Department of Mathematics, University of Illinois, Urbana-Champaign, Urbana, Illinois 61801, USA}
\email{heyizhu2@illinois.edu}

\date{April 2025}

\thanks{The first author was supported by NSF grant DMS-2202267.}



\begin{abstract}
We describe some periodic structure in the cohomology of the moduli stack of 1-dimensional formal group laws, also known as the $E_2$-page of the classical Adams--Novikov spectral sequence. This structure is distinct from the familiar $v_n$-periodicities, and it displays interesting number-theoretic properties. Our techniques involve the $\C$-motivic Adams spectral sequence, and we obtain analogous periodic structure in $\C$-motivic stable homotopy.
\end{abstract}

\maketitle
\tableofcontents

\section{Introduction}
\label{sctn:intro}

This manuscript exclusively considers $2$-primary stable homotopy.
In \cref{thm:C-homotopy-intro}, we describe some periodic families of non-trivial elements in $\C$-motivic stable homotopy.  Along the way in \cref{cor:Adams-Novikov-E2}, we describe the corresponding periodic families of non-trivial elements in the cohomology of the moduli stack of 1-dimensional formal group laws, also known as the classical Adams--Novikov $E_2$-page.

\subsection{The cohomology of the moduli stack of 1-dimensional formal group laws}

The moduli stack of 1-dimensional formal group laws is an object of interest in number theory.  The cohomology of this moduli stack is also known as the Adams--Novikov $E_2$-page, which is of interest to homotopy theorists because it is an approximation to stable homotopy via the Adams--Novikov spectral sequence.

The situation is well-illustrated by $v_1$-periodic homotopy.  In stems congruent to 7 modulo 8, $v_1$-periodic homotopy is detected in the Adams--Novikov $E_2$-page in filtration $1$.  (In stems congruent to 3 modulo 8, Adams--Novikov differentials and hidden extensions occur.)  In other words, $v_1$-periodic homotopy is detected by the cohomology of the moduli stack of 1-dimensional formal group laws.  This observation is intended to ``explain'' the appearance of number-theoretic structure in $v_1$-periodic homotopy.  Following this logic, we should not be surprised to find other regular patterns in the Adams--Novikov $E_2$-page that also display number-theoretic structure.  This manuscript describes such a regular pattern.

Long experience shows that the information in the stable homotopy groups of spheres can be sorted into separate periodic families.  It is a deep and central fact of classical stable homotopy theory that the only such $2$-primary periodicities are $v_n$ for $n \geq 0$.  The traditional $v_n$-periodicities appear along horizontal lines in the Adams--Novikov $E_2$-page. However, the Adams--Novikov $E_2$-page possesses additional algebraic periodicities \cite{Krause-thesis}.  We study a particular family of elements in the Adams--Novikov $E_2$-page that exhibits an ``exotic'' periodicity.  Our elements lie along lines of slope $\frac{1}{5}$ in the Adams--Novikov $E_2$-page.

\subsection{Motivation from topological modular forms}
\label{subsctn:mmf}

One of the many roles of the topological modular forms spectrum $\tmf$ in stable homotopy theory is that it serves as a testing ground for computational methods and strategies.  If a method computes the homotopy of $\tmf$ efficiently, then it is likely to be of use for the much more difficult sphere spectrum.

The manuscript \cite{IKLRZ} studies $\tmf$, and also its $\C$-motivic analogue $\mmf$, from a $g$-periodic perspective, where $g$ is a particular element in the 20-stem. The $g$-periodic structure of $\mmf$ reveals certain information that is otherwise difficult to obtain.  Therefore, we expect that a $g$-periodic perspective on the $\C$-motivic sphere is also interesting.

In this manuscript, building on previous work of Michael Andrews \cite{And}, we investigate $g$-periodic phenomena in the $\C$-motivic sphere spectrum.  We find a very rich structure, analogous to the number-theoretically interesting $v_1$-periodic homotopy, but quantifiably distinct from $v_1$-periodic phenomena.  In fact, our results only scratch the surface of a vastly more complicated structure.

\subsection{\texorpdfstring{$w_1$}{w1}-periodicity}

Remarkably, $\C$-motivic stable homotopy theory possesses additional periodicities that are completely separate from the classical $v_n$-periodicities \cite{Krause-thesis} \cite{MR3732063}.  A chromatic approach to $\C$-motivic homotopy theory necessarily requires that we study these exotic periodicities.

The first exotic periodicity is the Hopf map $\eta$, which is not nilpotent in $\C$-motivic homotopy theory.  The $\eta$-periodic $\C$-motivic stable homotopy groups are now well-understood, both computationally \cite{GI15} \cite{AM17} and theoretically \cite{BachmannHopkins}.  The second exotic periodicity is $w_1$, and it is the subject of this manuscript.  Note that $g$ from \cref{subsctn:mmf} is the same as $w_1^4$.

The goal of this manuscript is to describe a small portion of $\C$-motivic $w_1$-periodic homotopy that displays number-theoretically interesting properties.  Michael Andrews predicted that $\C$-motivic $w_1$-periodic homotopy would be at least as complicated as classical $v_1$-periodic homotopy.  Because $\C$-motivic $\eta$-periodicity is only slightly more complicated than classical $2$-periodicity, there is some hope that $w_1$-periodicity is only slightly more complicated than $v_1$-periodicity.  Our results show a very large increase in computational difficulty from the $v_1$-periodic context to the $w_1$-periodic context.  As of now, it is unrealistic to expect a complete description of the $w_1$-periodic homotopy of the $\C$-motivic sphere.

\subsection{Some \texorpdfstring{$w_1$}{w1}-periodic homotopy}

The $\C$-motivic homotopy groups display $w_1$-periodic structure that is analogous to the well-known classical $v_1$-periodic structure. We refer the reader ahead to \cref{subsctn:v1-homotopy} for an explicit recollection of the classical $v_1$-periodic computations.

The relationship of $v_1$ to $2$ is analogous to the relationship of $w_1$ to $\eta$.  For example, the classical cofiber of $2$ possesses a $v_1^4$ self-map, which leads to the familiar 8-fold periodicity of $v_1$-periodic homotopy.  Analogously, the $\C$-motivic cofiber of $\eta$ possesses a $w_1^4$ self-map of degree $(20,12)$ \cite{And}, which leads to infinite families of $\C$-motivic stable homotopy elements whose bidegrees are in arithmetic progression with difference $(20,12)$.  Here we grade $\C$-motivic stable homotopy groups by stem and motivic weight.

However, the $w_1^4$ self-map on $S/\eta$ is only useful for producing $w_1$-periodic families that are annihilated by $\eta$.  The most interesting parts of $v_1$-periodic homotopy consist of elements that are not annihilated by $2$.  Analogously, there is much to say about $\C$-motivic $w_1$-periodic homotopy elements that are not annihilated by $\eta$.  In this manuscript, we establish the existence of many $w_1$-periodic homotopy elements that are not annihilated by $\eta$.

In order to summarize our results, we will refer to the motivic coweight, which is the difference of the stem and the weight.  The coweight is relevant because $\eta$ has coweight zero.  Consequently, a homotopy element and its $\eta$ multiple share a common coweight, while their stems and weights both differ by $1$.

We describe non-zero $\C$-motivic stable homotopy elements in coweights that are congruent to 3 modulo 4.  We start with one previously identified element in each coweight that is annihilated by $\eta$, and we divide by $\eta$ multiple times to produce additional non-zero elements in the same coweight.  This is analogous to starting with the unique non-zero $2$-torsion element in a cyclic $v_1$-periodic stable homotopy group, and then dividing by $2$ repeatedly to determine the order of the cyclic group.

In many (but not all) coweights, we determine exactly how many times our $\eta$-torsion elements can be divided by $\eta$.  This gives us explicit submodules of the form $\Z[\eta]/(M, 2\eta, \eta^N)$ inside the $\C$-motivic stable homotopy groups $\pi_{*,*}$.  From another perspective, in many coweights, we produce non-zero elements that are not divisible by $\eta$, and we determine their $\eta$-exponents, i.e., the smallest value of $N$ such that $\eta^N$ annihilates the element.

Starting with coweight 3, we obtain the following $\eta$-exponents:
\begin{equation}
\label{eq:eta-exponents}
\infty, 3, 2, 7, 2, 6, 2, 15, 2, 6, 2, 12, 2, 6, 2, 31, 2, 6, 2, 12, 2, 6, 2, ?, 2, 6, 2, 12, 2, 6, 2, 63, \ldots.
\end{equation}
The first entry in the list is an exception.  In coweight 3, the element under consideration is $\sigma$ in $\pi_{7,4}$, which is $\eta$-periodic; no power of $\eta$ annihilates it.  In coweight $2^n - 1$, the element under consideration is Mahowald's $\eta_{n+1}$, and the $\eta$-exponent of $\eta_{n+1}$ is $2^{n-1}-1$.  In other words, the $2^{n-2}$th entry on the list equals $2^{n-1}-1$.

The remaining values display an interesting structure.  Every other entry equals $2$; every fourth entry equals $6$; and every eighth entry equals $12$.  In other words, except for the special cases described in the above paragraph, the $k$th value on the list depends only on the $2$-adic valuation of $k$ in many cases.  We expect that this number-theoretic behavior is a general phenomenon, although we have not proved it in general.

The first unknown $\eta$-exponent occurs at the $24$th entry on the list, i.e., in coweight $95 = 4 \cdot 24 - 1$.  We expect that the same (currently unknown) value occurs in every sixteenth entry.

The previous discussion narrates our main results.  We now state these results in the form of a theorem, although the careful reader should inspect \cref{thm:C-homotopy} and \cref{thm:C-homotopy-beta} for more precise statements.

\begin{thm} 
(see \cref{thm:C-homotopy} and \cref{thm:C-homotopy-beta})
\label{thm:C-homotopy-intro}
For all $n \geq 1$, there is a non-zero element in $\C$-motivic stable homotopy such that:
\begin{enumerate}
\item
the coweight of the element is $4n-1$.
\item 
the element is not divisible by $\eta$.
\item
when $n = 2^j$, the element equals (the $\C$-motivic version of) the Mahowald element $\eta_{j+3}$
in stem $2^{j+3}$.
\item 
for $n \geq 2$, 
the $n$th term of sequence \eqref{eq:eta-exponents} gives the smallest number $N$ such that $\eta^N$ annihilates the element.
\end{enumerate}
\end{thm}

\begin{rem}
\cref{thm:C-homotopy-intro} can be compared to the following well-known claims about classical $v_1$-periodicity.  For all $n \geq 1$, there is a non-zero element in classical stable homotopy such that:
\begin{enumerate}
\item
the stem of the element is $4n-1$.
\item 
the element is not divisible by $2$.
\item 
the $n$th term of the sequence
\[
3, 4, 3, 5, 3, 4, 3, 6, 3, 4, 3, 5, 3, 4, 3, 7, \ldots
\]
gives the smallest number $N$ such that $2^N$ annihilates the element.
\end{enumerate}
See \cref{subsctn:v1-homotopy} below for a further description of this analogous classical $v_1$-periodic structure.
\end{rem}

\cref{fig:w1-Adams-Einfty} illustrates our main theorems in the form of an Adams chart.  The chart shows the elements described in \cref{thm:C-homotopy-intro} in red.  The elements in gray are $w_1$-periodic elements previously identified by Andrews \cite{And}.

The key steps in establishing part (4) of \cref{thm:C-homotopy-intro} are specific Adams differential computations.  This information affects the values of the terms in sequence \eqref{eq:eta-exponents}.  First, we find an infinite $w_1$-periodic family of non-zero differentials
\[
d_2( h_3 g^k) = h_0 h_2 \cdot h_2 g^k.
\]
In order to obtain this family of differentials, we must do some extensive algebraic preparatory work; see \cref{subsctn:BX-intro} for more discussion.  Second, we also establish a specific permanent cycle in the $\C$-motivic 107-stem.  This requires detailed analysis of elements in nearby stems; see \cref{subsctn:D1h12i1} for more details.

\begin{rem}
As an aside, we point out that the classical Adams spectral sequence for $\tmf$ can be used to detect additional infinite $w_1$-periodic families of non-zero differentials.  For example, the element $d_0 g^4$ in the 94-stem is the value of a classical Adams $d_2$ differential \cite{IWX22b}.  Comparison to the Adams spectral sequence for $\tmf$ shows that $d_0 g^k$ is non-zero in the Adams $E_2$-page for all $k$.  Therefore, the elements $d_0 g^k$, for $k \geq 4$, are the targets of a $w_1$-periodic family of classical differentials.

In the $\C$-motivic context, there is a corresponding family of Adams $d_2$ differentials whose targets are $\tau^6 d_0 g^k$ for $k \geq 4$.
\end{rem}

\subsection{The Adams--Novikov \texorpdfstring{$E_2$}{E2}-page}

This manuscript describes a family of non-zero elements in the Adams--Novikov $E_2$-page, whose degrees increase arithmetically with common difference $(20,4)$.  Each of these elements is annihilated by the element $\alpha_1$ that detects $\eta$.  We divide by $\alpha_1$ multiple times to produce additional non-zero elements.  In many (but not all) cases, we determine exactly how many times our $\alpha_1$-torsion elements are divisible by $\alpha_1$.  From another perspective, we produce non-zero elements that are not divisible by $\alpha_1$, and we determine their $\alpha_1$-exponents, i.e., the smallest value of $N$ such that $\alpha_1^N$ annihilates the element.  These exponents are displayed in the sequence \eqref{eq:eta-exponents}.

All of our results about the classical Adams--Novikov $E_2$-page are corollaries of our study of $\C$-motivic homotopy.  See \cite[Chapter 6]{Isa19} for the close connection between the classical Adams--Novikov $E_2$-page and $\C$-motivic homotopy.  In particular, \cite[Theorem 6.12]{Isa19} equates the Adams--Novikov $E_2$-page with the homotopy of the $\C$-motivic two-cell complex $S/\tau$.  Our study of $\C$-motivic homotopy yields elements that are all detected by the inclusion $S \rightarrow S/\tau$.  So we learn about corresponding elements in $S/\tau$, and from there to the Adams--Novikov $E_2$-page.  See \cref{sctn:Adams-Novikov} for additional details of the translation between $\C$-motivic homotopy and the classical Adams--Novikov $E_2$-page.  We state the translation of \cref{thm:C-homotopy-intro} here, although the reader should inspect \cref{thm:Adams-Novikov-E2} and \cref{thm:Adams-Novikov-E2-beta} for more precise statements.

\begin{landscape}
\begin{figure}[H]
\begin{center}
\includegraphics[trim={0cm, 0cm, 0cm, 0cm},clip,page=1,scale=0.77]{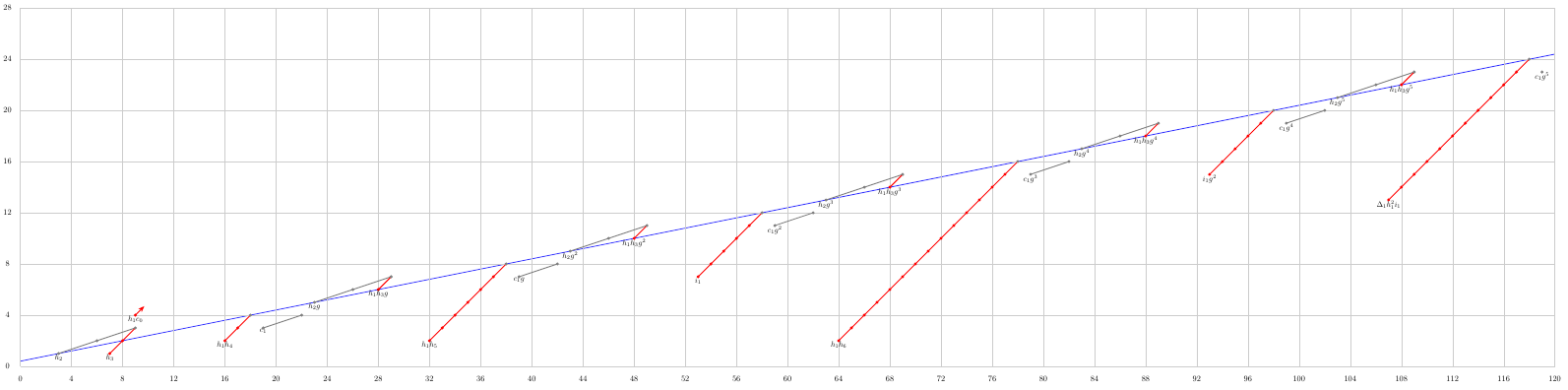}
\caption{Some $w_1$-periodic elements in the $\C$-motivic Adams $E_\infty$-page.  Gray elements were previously studied by Andrews.  Red elements are studied in this manuscript.
\label{fig:w1-Adams-Einfty}}
\end{center}
\end{figure}
\end{landscape}

\begin{cor}
(see \cref{thm:Adams-Novikov-E2} and \cref{thm:Adams-Novikov-E2-beta})
\label{cor:Adams-Novikov-E2}
For all $n \geq 1$, there is a non-zero element in the cohomology of the moduli stack of 1-dimensional formal group laws such that:
\begin{enumerate}
\item 
the difference between the stem and the cohomological filtration of the element is $8n-2$.
\item
the element is not divisible by $\alpha_1$.
\item 
When $n = 2^j$, the element equals the  $\beta$-family element that detects Mahowald's $\eta_{j+3}$ in stem $2^{j+3}$.
\item 
For $n \geq 2$, the $n$th term of sequence \eqref{eq:eta-exponents} gives the smallest number $N$ such that $\alpha_1^N$ annihilates the element.
\end{enumerate}
\end{cor}

\cref{fig:w1-Adams-Einfty} demonstrates the contents of \cref{thm:C-homotopy-intro}, but it also serves as an illustration of \cref{cor:Adams-Novikov-E2}.  Almost all of the elements that we study, as well as Andrews's families, share the property that their Adams filtrations equal their Adams--Novikov filtrations, so the Adams chart of \cref{fig:w1-Adams-Einfty} is also an Adams--Novikov chart.  However, the Adams filtrations and Adams--Novikov filtrations differ in the 9-stem.  In the case of the Adams--Novikov $E_2$-page, there are two distinct non-zero elements in filtration 3.  One element is a multiple of $\alpha_{2/2}$.  The other element is a multiple of $\alpha_1^2$, and it supports infinitely many $\alpha_1$-multiplications.

\subsection{Relationship to Mahowald's \texorpdfstring{$\eta_j$}{eta j} elements}
The $\C$-motivic versions of Mahowald's $\eta_j$ elements appear as distinguished elements amongst the $w_1$-periodic homotopy that we study.  We find that 
\begin{equation}
\label{eq:etaj-family}
w_1^{2^{j-2}} \cdot \eta_j = \eta^{2^{j-2}} \eta_{j+1}.
\end{equation}
In other words, each $\eta_j$ is the first element of a $w_1^{2^{j-2}}$-periodic family, and the second element of this family is an $\eta$-multiple of the next $\eta_{j+1}$.  We see here a ``family'' structure on the Mahowald $\eta_j$'s that we believe has not been previously observed.

The structure on the $\eta_j$ family is analogous to the structure exhibited by the classical elements $\rho_{2^n-1}$ that generate $v_1$-periodic homotopy in the $(2^n-1)$-stem.  We have 
\[
v_1^{2^{n-1}} \cdot \rho_{2^n-1} = 2 \rho_{2^{n+1}-1}.
\]
Note the factor of $2$ on the right side of the formula, corresponding to the powers of $\eta$ on the right side of \cref{eq:etaj-family}.

\subsection{Relationship to Andrews's work}

Michael Andrews \cite{And} previously studied several infinite families of $w_1$-periodic elements.  Of the elements discussed in \cref{thm:C-homotopy-intro}, only the $\eta$-torsion elements were found by Andrews (i.e., one element in each coweight).  Consequently, Andrews's work does not exhibit number-theoretic properties.

Parts (1), (2), and (3) of our main Theorems \ref{thm:C-homotopy} and \ref{thm:C-homotopy-beta} were first established by Andrews \cite{And}, although we reprove these parts independently.  On the other hand, part (4) of both theorems is genuinely new and requires detailed, technical work.

\begin{figure}[H]
\begin{center}
\includegraphics[trim={0cm, 0cm, 0cm, 0cm},clip,page=1,scale=1]{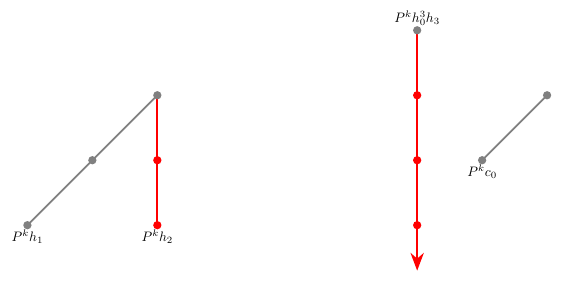}
\caption{The classical $v_1$-periodic motif of period $(8,4)$.  Gray dots detect $2$-torsion; red dots detect elements that are not $2$-torsion.
\label{fig:v1-motif}}
\end{center}
\end{figure}

\cref{fig:v1-motif} shows one period of the familiar $v_1$-periodic region at the top of an Adams chart.  The arrow indicates that the columns have varying lengths in the $(8k-1)$-stem.  Note also that the gray dots are precisely the $h_0$-torsion elements.

\begin{figure}[H]
\begin{center}
\includegraphics[trim={0cm, 0cm, 0cm, 0cm},clip,page=1,scale=0.64]{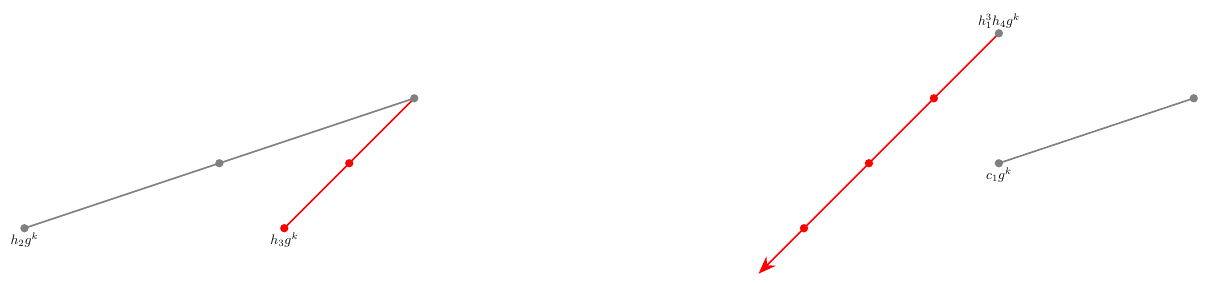}
\caption{A $\C$-motivic $w_1$-periodic motif of period $(20,4,12)$.  Gray dots detect $h_1$-torsion; red dots detect elements that are not $h_1$-torsion.
\label{fig:w1-motif}}
\end{center}
\end{figure}

The analogy between $v_1$-periodicity and $w_1$-periodicity involves a degree shift in terms of Adams charts.  Consequently, the classical $v_1$-periodic motif of \cref{fig:v1-motif} corresponds to the skewed $w_1$-periodic motif shown in \cref{fig:w1-motif}.  Andrews studied the families shown in gray; we study the families shown in red.

\subsection{The Burklund--Xu spectral sequence}
\label{subsctn:BX-intro}

Our main algebraic tool is the Burk\-lund--Xu spectral sequence \cite{BX23}. This spectral sequence was originally intended to compute the $E_2$-page of the Cartan--Eilenberg spectral sequence associated to the  extension
\[
\F_2[\zeta_1^2, \zeta_2^2, \ldots] \rightarrow
\F_2[\zeta_1, \zeta_2, \ldots] \rightarrow 
E(\zeta_1, \zeta_2, \ldots),
\]
where the middle object is the classical dual Steenrod algebra, and the quotient object on the right is an exterior algebra.  In order to obtain information about the Adams $E_2$-page, one must additionally study the Cartan--Eilenberg differentials.

Later we will be particularly interested in elements of Chow degree one, i.e., elements in degree $(s,f,w)$ such that $s+f-2w=1$.  Here $s$, $f$, and $w$ are the stem, Adams filtration, and motivic weight respectively.

It turns out that for elementary degree reasons, there are no possible $\C$-motivic Cartan--Eilenberg differentials in Chow degree one.  Consequently, the $\C$-motivic Burk\-lund--Xu spectral sequence converges all the way to the Adams $E_2$-page in Chow degree one.

The Burklund--Xu spectral sequence is a remarkably powerful tool for studying the $\C$-motivic Adams $E_2$-page in Chow degree one.  One may reasonably ask about the significance of computations in Chow degree one.  In our case, there is a particular family of products $h_0 h_2 \cdot h_2 g^k$ in degrees $(20k+6, 4k+3, 12k+4)$ of relevance to our study of $w_1$-periodic homotopy.  We use the Burklund--Xu spectral sequence to show that all of these products are non-zero.

We draw attention to \cref{fig:A-E2} and \cref{fig:A-Einfty}, which display the Burklund--Xu spectral sequence in Chow degree one in a range, including all differentials and some hidden extensions in its $E_\infty$-page.  These charts demonstrate that the spectral sequence is computable in practice.

The related article \cite{KILRZ-e0gk} offers another demonstration of the power of the Burklund--Xu spectral sequence.  It studies another infinite family of elements in Chow degree one.

\subsection{Future directions}

Our work suggests a number of directions for further study.

\begin{problem}
\label{prob:unknown-exponent}
Determine the first unknown value in the sequence \eqref{eq:eta-exponents}.  In other words, compute the $\eta$-exponent of the element in coweight $95 = 4 \cdot 24 - 1$.
\end{problem}

We do not have a practical plan of attack for \cref{prob:unknown-exponent}.  We handle the previous case in coweight $47 = 4 \cdot 12 - 1$ by a detailed analysis of the $\C$-motivic Adams spectral sequence around stem $107$.  That approach is unlikely to generalize, so a genuinely new idea may be required.

\begin{problem}
\label{prob:unknown-exponent-general}
What is the $\eta$-exponent in coweight $3 \cdot 2^n - 1$?     
\end{problem}

\begin{problem}
\label{prob:adic}
Prove that the values of $\eta$-exponents in coweight $4k-1$ depend only on the $2$-adic valuation of $k$, with the exception of coweights of the form $2^n-1$.
\end{problem}

This manuscript presents some evidence to suggest that the exponents do in fact depend only on the 2-adic valuations.  Solutions to both \cref{prob:unknown-exponent-general} and \cref{prob:adic} would completely determine the entire sequence \eqref{eq:eta-exponents}.

\begin{problem}
Give a number-theoretic description of the sequence \eqref{eq:eta-exponents}.
\end{problem}

Sequence \eqref{eq:eta-exponents} describes structure in the cohomology of the moduli stack of 1-dim\-en\-sion\-al formal group laws, so it is plausible that such a number-theoretic description exists.

\begin{problem}
Compute more $\C$-motivic $w_1$-periodic homotopy.
\end{problem}

Our work just scratches the surface.  Naive inspection turns up additional $w_1$-periodic elements, but we have not organized these examples into any kind of coherent structure.

\begin{problem}
\label{prob:H1}
Consider the element $H_1$ in degree $(62,5,33)$.  Determine the $w_1$-periodic behavior of this element.
\end{problem}

The element $H_1$ appears in \cref{fig:Guozhen} at location $(29,5)$.  The figure suggests that $H_1$ is $w_1^4$-periodic in the Adams $E_2$-page, but we have not established this guess.  We guess that the elements $h_2 H_1 g^k$, $h_3 H_1 g^k$, $D_4$, and $h_2 G_0$ are similarly periodic.  We expect that all of these elements are of interest in a $w_1$-periodic perspective, but we do not know why, nor how to predict additional interesting elements.

\begin{problem}
Is there a $w_1$-periodic version of the classical image of $J$ spectrum?
\end{problem}

Such a ``$\mathit{wj}$'' spectrum would cleanly detect the $w_1$-periodic homotopy that we explore in this manuscript.  One could start with the $\mathit{wko}$ spectrum \cite{Quigley} whose homotopy is
\[
\frac{\F_2[\eta, \nu, u, \kappabar]}{\eta \nu, \nu^3, \nu u, u^2 + \eta^2 \kappabar},
\]
then construct a self-map of $\mathit{wko}$, analogous to the classical $1 - \psi^3$, whose fiber has the desired properties for $\mathit{wj}$.

\begin{problem}
\label{prob:BXss}
Determine the differentials and hidden extensions in the Burklund--Xu spectral sequence in Chow degree one in a range, thereby determining the Adams $E_2$-page in Chow degree one in that range.
\end{problem}

We suspect that  \cref{prob:BXss} is manageable at least as far as the range displayed in \cref{fig:Guozhen}, although it would require substantial effort to do so by hand.  It is likely that machines could carry out the computation more effectively.

\subsection{Outline}

\cref{sctn:background} provides brief recollections about our degree conventions, the structure of the dual Steenrod algebra, and some essential results about the cohomology of the Steenrod algebra.  It also introduces the notation and terminology to be used later.

\cref{sctn:v1-w1} describes some $v_1$-periodic structure in the cohomology of the classical Steenrod algebra.  We study families of elements in stems congruent to 3 modulo 8 and 7 modulo 8.  Here we are studying algebraic phenomena, in the sense that we consider structure in the Adams $E_2$-page rather than the Adams $E_\infty$-page.  \cref{sctn:v1-w1} also describes analogous algebraic $w_1$-periodic structure in the cohomology of the $\C$-motivic Steenrod algebra, i.e., in the $\C$-motivic Adams $E_2$-page.

\cref{sctn:v1-w1-periodic-homotopy} starts with a brief summary of the well-known $v_1$-periodic homotopy groups in \cref{subsctn:v1-homotopy}, and we describe how they are detected in the Adams $E_\infty$-page in \cref{subsctn:v1-periodic-homotopy-Adams}.  Then \cref{subsctn:w1-periodic-homotopy} contains the statements and proofs of our main results.  The previous section described structure in the $\C$-motivic Adams $E_2$-page, and this section describes the fate of this $w_1$-periodic structure in the Adams spectral sequence.  The proofs of the most difficult and interesting parts of our main theorems are deferred until later because they require extensive technical work with specific elements.

\cref{sctn:BXss} begins the technical work required to prove the hard parts of the theorems in \cref{subsctn:w1-periodic-homotopy}.  We discuss the $\C$-motivic Burklund--Xu spectral sequence, which gives us firm control over the $\C$-motivic Adams $E_2$-page in Chow degree one.

\cref{sctn:BXss-analysis} describes our technical work with the $\C$-motivic Burklund--Xu spectral sequence.  We study specific classes in specific degrees, leading to the existence of an infinite family of non-zero products $h_0 h_2 \cdot h_2 g^k$ in the $\C$-motivic Adams $E_2$-page.  This family is relevant later in our analysis of Adams differentials.

\cref{sctn:Adams-differentials} calculates some $\C$-motivic Adams differentials.  We find one infinite family 
\[
d_2(h_3 g^k) = h_0 h_2 \cdot h_2 g^k
\]
of non-zero differentials.  We also establish several infinite families of non-zero permanent cycles.  These results about Adams differentials allow us to finish the proofs of our main theorems from \cref{subsctn:w1-periodic-homotopy}.

\cref{sctn:Adams-Novikov} gives a restatement of our main results in different language.  There is a straightforward translation between the structure of $\C$-motivic homotopy and the structure of the cohomology of the moduli stack of 1-dimensional formal group laws.  This section simply translates the theorems of \cref{subsctn:w1-periodic-homotopy}.

\cref{sctn:A(2)} is a digression from the main story.  We carry out an exhaustive study of the part of the cohomology of $\C$-motivic $\cA(2)$ in Chow degree one.  Our digression illustrates a general philosophy that one should always study the cohomology of $\cA(2)$ as a warmup for the vastly more difficult cohomology of the Steenrod algebra.  This principle applies in many contexts, including our specific interest in Chow degree one in this manuscript.  This section can serve as a training module for readers who seek hands-on experience with the Burklund--Xu spectral sequence. 

\section{Background}
\label{sctn:background}

Regarding gradings,
\begin{itemize}
\item
in motivic bigraded contexts, such as in $\C$-motivic stable homotopy, all elements are graded in the form $(s,w)$, where $s$ is the topological stem and $w$ is the motivic weight.
\item
in motivic trigraded contexts, such as in the $\C$-motivic Adams $E_2$-page, all elements are graded in the form $(s,f,w)$, where $f$ is the Adams filtration.
\item
in classical bigraded contexts, such as the classical Adams $E_2$-page, all elements are graded in the form $(s,f)$.
\item
the Chow degree of a trigraded element is $s+f-2w$, and the Chow degree of a motivic bigraded element is $s-2w$.
\item 
we also use the motivic coweight, which is defined to be $s-w$, in specific situations where that grading is convenient.
\item
for an element in stem $s$ and Adams filtration $f$, the $v_1$-intercept is the number $s-2f$.  On a standard Adams chart, the $v_1$-intercept is the $x$-intercept of the line of slope $\frac{1}{2}$ that passes through the point $(s,f)$.  The use of $v_1$ in the terminology reflects that the slope of $v_1$-multiplication is $\frac{1}{2}$.
\item
\cref{sctn:BX-Cd1} discusses an ad hoc grading system that we use for the Burklund--Xu spectral sequence in Chow degree one.
\end{itemize}

The $\C$-motivic dual Steenrod algebra $\cA_{**}$ is 
\[
\frac{\mathbb{F}_2 [\tau][\tau_i,\xi_j | i\geq 0, j\geq 1]}{\tau_i^2=\tau\xi_{i+1}},
\]
where the internal degrees and weights of $\tau_i$ and $\xi_j$ are $(2^{i+1}-1,2^i-1)$ and $(2^{j+1}-2,2^j-1)$ respectively \cite{voev03}.  The ground ring of this Hopf algebra is $\F_2[\tau]$, which is the $\C$-motivic homology of a point.

The classical dual Steenrod algebra $\cA^{\cl}_*$ is
\[
\mathbb{F}_2[\zeta_i | i\geq 1],
\]
where $\zeta_i$ has degree $2^i - 1$.  Betti realization (also known as $\tau$-localization) takes $\tau_i$ to $\zeta_{i+1}$ and takes $\xi_j$ to $\zeta_j^2$.

We use $\cA(1)_{**}$ and $\cA(2)_{**}$ for the usual quotients of $\cA_{**}$, and we use $\cA^{\cl}(1)_*$ and $\cA^{\cl}(2)_*$ for the usual quotients of $\cA^{\cl}_*$.

We write $H^{***} \cA$ for the trigraded cohomology of the $\C$-motivic Steenrod algebra $\cA$.  This cohomology is also known as the $\C$-motivic Adams $E_2$-page, and it can be described homologically as $\Ext_{\cA}(\F_2[\tau], \F_2[\tau])$.  We write $H^{***} \cA(1)$ and $H^{***} \cA(2)$ for the trigraded cohomologies of $\cA(1)$ and $\cA(2)$.

We write $H^{**} \cA^{\cl}$ for the bigraded cohomology of the classical Steenrod algebra $\cA^{\cl}$.  This cohomology is also known as the classical Adams $E_2$-page, and it can be described homologically as $\Ext_{\cA^{\cl}}(\F_2, \F_2)$.  We write $H^{**} \cA^{\cl}(1)$ and $H^{**} \cA^{\cl}(2)$ for the bigraded cohomologies of $\cA^{\cl}(1)$ and $\cA^{\cl}(2)$.

We write $\cB_{**}$ for the subalgebra $\F_2[\xi_j | j \geq 1]$ of $\cA_{**}$.  Beware that the ground ring $\F_2$ of $\cB_{**}$ is not equal to the ground
ring $\F_2[\tau]$ of $\cA_{**}$.  The most important point about $\cB_{**}$ is that its cobar complex is precisely equal to the part of the cobar complex of $\cA_{**}$ in Chow degree zero.  Moreover, there is an isomorphism $\cA^{\cl}_* \rightarrow \cB_{**}$ of Hopf algebras that takes $\zeta_i$ to $\xi_i$.  Beware that the isomorphism takes elements of internal degree $t$ to elements of internal degree $2t$ and motivic weight $t$.

We write $H^{***} \cB$ for the trigraded cohomology of the dual $\cB$ of $\cB_{**}$.  However, $H^{***} \cB$ is concentrated in Chow degree zero, i.e.,
$H^{s,f,w} \cB = 0$ if $s+f-2w \neq 0$.  This follows immediately from the observation that the cobar complex of $\cB$ is concentrated in Chow degree zero.

\begin{thm}
\cite[Theorem 2.19]{Isa19}\label{Ddi}
There are highly structured isomorphisms be\-tween $H^{**} \cA^{\cl}$, $H^{***} \cB$, and the part of $H^{***} \cA$ that lies in Chow degree zero.
\end{thm}

An explicit version of \cref{Ddi} is that
\begin{equation}
\label{eq:Ddi}
H^{s,f,w}\cB \cong 
\begin{cases}
H^{\frac{s-f}{2},f}\cA^{\cl} & \text{if}~s+f-2w=0.\\
0 &\text{otherwise.}
\end{cases}
\end{equation}
\cref{eq:Ddi} is particularly useful for keeping track of the degree shifts in the isomorphism. 

We will offer another proof of \cref{Ddi} in \cref{proof:Cd0} that uses the techniques that are developed in this manuscript.  The term ``highly structured isomorphism''  is somewhat vague.  The point is that all three objects are the homologies of isomorphic differential graded algebras, so their Yoneda products, squaring operations, and Massey products are all identical.

We write $1^{BX} \cdot (-)$ for the isomorphism from $H^{**} \cA^{\cl}$ to the part of $H^{***} \cA$ that lies in Chow degree zero.  This seemingly obscure notation makes more sense in the context of the $\C$-motivic Burklund--Xu spectral sequence to be discussed below in \cref{sctn:BXss}.

\begin{prop}
\label{prop:Chow-degree-nonneg}
Every non-zero element of $H^{***}\cA$ has a non-negative Chow degree.
\end{prop}

\begin{proof}
By inspection of degrees, the cobar complex for $\cA$ is concentrated in non-negative Chow degrees.
\end{proof}

\subsection{Submodules, cogenerators, and coexponents}
\label{subsctn:coexponent}

Let $a$ be an element of a (typically graded) commutative $\Z$-algebra $A$.  To avoid pathologies, we assume that no non-zero element of $A$ is infinitely divisible by $a$.  In practice, this kind of condition is typically guaranteed by a bounded-below hypothesis on the graded algebra $A$.

An $a$-submodule of $A$ is $\Z[a]$-submodule of $A$, i.e., a subgroup of the underlying abelian group of $A$ that is closed under $a$-multiplication.  We are primarily interested in cyclic $a$-submodules, i.e., $a$-submodules that are generated by a single element.

In general, there are many isomorphism types of cyclic $a$-submodules, such as
\[
\frac{\Z[a]}{2^{11}, 2^6 a, 2^5 a^2}.
\]
However, in practice we will only consider elements $a$ such that $2 a = 0$.  With this additional condition on $a$, every cyclic submodule is of the form
\[
\frac{\Z[a]}{M, 2 a, a^N}
\]
for some value of $M$ in $\N$ and some value of $N$ in $\N \cup \{\infty\}$.  The value of $M$ is the characteristic of the $a$-submodule, but we will generally not study $M$.  We refer to $N$ as the \emph{$a$-exponent} of the cyclic $a$-submodule.

Let $\frac{x}{a}$ be an $a$-torsion element of $A$.  The notation $x$ itself (without a denominator) is not an element of $A$.  Rather, it refers to the relation
\[
a \cdot \frac{x}{a} = 0.
\]
This naming convention may seem confusing at first, but it has structural advantages in our applications.  It is similar to the notation that is typically used to describe $v_1$-periodic elements of the Adams--Novikov $E_2$-page in filtration 1 \cite[Chapter 5.2]{Ravenel86}.  This similarity is not coincidental since we also describe $v_1$-periodic phenomena.

More generally, we write $\frac{x}{a^n}$ for any element $y$ with the property that $a^{n-1} \cdot y = \frac{x}{a}$.  Beware that the element $\frac{x}{a^n}$ may not exist if $n$ is too large for fixed $x$, and it is not necessarily unique if there happen to be $a^{n-1}$-torsion elements.  This ambiguity is no problem for us.  In practice, our arguments apply equally well to any choice of $\frac{x}{a^n}$.

The element $\frac{x}{a}$ is the unique non-zero $a$-torsion element inside of a largest cyclic $a$-submodule $C$ of $A$, where $C$ is isomorphic to $\Z[a]/(M, 2a, a^N)$ for some $M$ and $N$.  We refer to $\frac{x}{a}$ as the \emph{cogenerator} of $C$, and we refer to $N$ as the \emph{$a$-coexponent} of $\frac{x}{a}$.  The $a$-coexponent equals one more than the number of times that $a$ divides the element $\frac{x}{a}$.

\begin{ex}
In \cref{subsctn:v1-periodic-h0-submodules}, we study $h_0$-submodules of the classical Adams $E_2$-page.  Here we are working in an $\F_2$-algebra, so the characteristic of every $h_0$-submodule is $2$.
\end{ex}

\begin{ex}
In \cref{subsctn:w1-periodic-h1-submodules}, we study $h_1$-submodules of the $\C$-motivic Adams $E_2$-page.  Again we are working in an $\F_2$-algebra, so the characteristic of every $h_1$-submodule is $2$.
\end{ex}

\begin{ex}
\label{ex:eta-submodules}
In \cref{subsctn:w1-periodic-homotopy}, we study $\eta$-submodules of the $\C$-motivic stable homotopy ring $\pi_{*,*}$.  Now we are not working in an $\F_2$-algebra, so the characteristic can be larger than $2$.  In practice, some of the $\eta$-submodules that we study do have characteristics larger than $2$, but we will not study characteristics because we do not have techniques to analyze them.  Rather, we are more interested in $\eta$-exponents.  Note that $2 \eta = 0$, so our simpler classification of $\eta$-submodules applies.
\end{ex}

\begin{ex}
\label{ex:alpha1-submodules}
In \cref{sctn:Adams-Novikov}, we study $\alpha_1$-submodules of the cohomology of the moduli stack of 1-dimensional formal group laws, also known as the classical Adams--Novikov $E_2$-page.  Again, we are not working in an $\F_2$-algebra, so the characteristic can be larger than $2$, but we will not study characteristics.  Rather, we are more interested in $\alpha_1$-exponents.  Note that $2 \alpha_1 = 0$, so our simpler classification of $\alpha_1$-submodules applies.
\end{ex}

\section{Classical \texorpdfstring{$v_1$}{v1}-periodicity and \texorpdfstring{$\C$}{C}-motivic \texorpdfstring{$w_1$}{w1}-periodicity in algebra}
\label{sctn:v1-w1}

\subsection{\texorpdfstring{$v_1$}{v1}-periodicity in the cohomology of the classical Steenrod algebra}
\label{subsctn:v1-algebraic}

The classical $v_1$-periodic stable homotopy groups are represented in the classical Adams spectral sequence by specific elements in high Adams filtration, near the vanishing line of slope $\frac{1}{2}$.  The goal of this section is to give a concrete description of these representatives in the Adams spectral sequence.  This material is well-known.  We will generalize all of this structure into the more exotic $\C$-motivic context, so we need a thorough understanding of the classical situation first.

Let $\frac{\overline{v_1^{4k}}}{h_0}$ be the unique non-zero element in $H^{**}\cA^\mathrm{cl}$ in degree $(8k-1, 4k)$.  This element is easily identified on a standard Adams chart such as \cite{IWX22b}.  It is the element at the top of the $h_0$-tower, at the top of the $(8k-1)$-stem.  The notation reflects the fact that $v_1^{4k}$ is the reason that $\frac{\overline{v_1^{4k}}}{h_0}$ is annihilated by $h_0$; in other words, in the various spectral sequences that compute $H^{**}\cA^{\cl}$, there is an element, typically called $v_1^{4k}$ among other possible names, that supports a differential whose value is $h_0$ times a non-zero element in degree $(8k-1,4k)$.

We use the terminology and notation of \cref{subsctn:coexponent} to describe $h_0$-submodules in terms of cogenerators and $h_0$-coexponents, so $\frac{\overline{v_1^{4k}}}{h_0^n}$ is an element $y$ with the property that $h_0^{n-1} \cdot y$ equals $\frac{\overline{v_1^{4k}}}{h_0}$.

\begin{rem}
\label{rem:v1-brackets}
The elements $\frac{\overline{v_1^{4k}}}{h_0}$ are closely related to each other via the usual Massey products that express $v_1$-periodicity.  Specifically, we have
\[
\frac{\overline{v_1^{4j+4k}}}{h_0} = 
\left\langle \frac{\overline{v_1^{4j}}}{h_0}, h_0, \frac{\overline{v_1^{4k}}}{h_0} \right\rangle
\]
with no indeterminacy for all $j$ and $k$.

More generally, assume that $n \leq m$.  For each choice of $\frac{\overline{v_1^{4j}}}{h_0^m}$ and each choice of $\frac{\overline{v_1^{4k}}}{h_0^n}$, the Massey product
\[
\left\langle \frac{\overline{v_1^{4j}}}{h_0^m}, h_0^m, \frac{\overline{v_1^{4k}}}{h_0^n} \right\rangle
\]
has no indeterminacy and equals an element that is a possible choice for $\frac{\overline{v_1^{4j+4k}}}{h_0^n}$.  Beware that the unique element in the Massey product may depend on the choices of $\frac{\overline{v_1^{4j}}}{h_0^m}$ and of $\frac{\overline{v_1^{4k}}}{h_0^n}$.

For a fixed value of $n$, consider the elements of the form $\frac{\overline{v_1^{4k}}}{h_0^n}$ in $H^{**}\cA^{\cl}$.  Such an element exists only when the $h_0$-coexponent of $\frac{\overline{v_1^{4k}}}{h_0}$ is at least $n$.  This family of existing elements forms a $v_1$-periodic family, but the degree of the periodicity depends on $n$.
\end{rem}

\begin{lem}
\label{lem:h0-coexponent-2^j}
For all $j \geq 2$, the $h_0$-coexponent of $\frac{\overline{v_1^{2^j}}}{h_0}$ is $2^{j-1}$.
\end{lem}

\begin{proof}
This follows immediately from the fact that
\[
h_0^{2^{j-1}-1} \cdot h_{j+1} = \frac{\overline{v_1^{2^j}}}{h_0}.
\]
\end{proof}

\begin{lem}
\label{lem:h0-coexponents}
For all $k \geq 0$, the $h_0$-coexponent of:
\begin{enumerate}
\item
$\frac{\overline{v_1^{8k+12}}}{h_0}$ is $6$.
\item 
$\frac{\overline{v_1^{16k+24}}}{h_0}$ is $12$.
\item 
$\frac{\overline{v_1^{32k+48}}}{h_0}$ is $28$.
\end{enumerate}
\end{lem}

\begin{proof}
In all three cases, the proof follows a similar logic.  Determine the $h_0$-co\-ex\-pon\-ents in low dimensions by explicit computation.  Then Adams periodicity \cite{May-never} \cite{Li20} determines the $h_0$-co\-ex\-pon\-ents in higher dimensions.

For part (1), the $h_0$-coexponents of $\frac{\overline{v_1^{12}}}{h_0}$, $\frac{\overline{v_1^{20}}}{h_0}$, and $\frac{\overline{v_1^{28}}}{h_0}$ are 6 by explicit computation \cite{IWX22b}.  For $k \geq 3$, we use that the Adams periodicity operator $P^4$ of degree $(32,16)$ induces an isomorphism on $H^{**}\cA^\mathrm{cl}$ in the relevant range.

For part (2), the $h_0$-coexponent of $\frac{\overline{v_1^{24}}}{h_0}$ is 12 by explicit computation \cite{IWX22b}.  For $k \geq 1$, we use that the Adams periodicity operator $P^4$ of degree $(32,16)$ induces an isomorphism on $H^{**}\cA^\mathrm{cl}$ in the relevant range.

For part (3), the $h_0$-coexponents of $\frac{\overline{v_1^{48}}}{h_0}$ and $\frac{\overline{v_1^{80}}}{h_0}$ are 28 by explicit computation \cite{IWX22b} \cite{BFCC23}.  For $k \geq 2$, we use that the Adams periodicity operator $P^8$ of degree $(64,32)$ induces an isomorphism on $H^{**}\cA^\mathrm{cl}$ in the relevant range.
\end{proof}

\begin{lem}
\label{lem:h0-coexponent-96}
The $h_0$-coexponent of $\frac{\overline{v_1^{96}}}{h_0}$ is 58.
\end{lem}

\begin{proof}
This specific case is explicitly computed in \cite{BFCC23}.
\end{proof}

\begin{rem}
\label{rem:h0-coexponent-96}
Unlike \cref{lem:h0-coexponents}, \cref{lem:h0-coexponent-96} only considers a single coexponent rather than a periodic family.  The problem is that the generator $\frac{\overline{v_1^{96}}}{h_0^{58}}$ lies slightly outside of the range in which (some power of) Adams periodicity is known to be an isomorphism.  If we could verify that the $h_0$-coexponent of $\frac{\overline{v_1^{64+96}}}{h_0}$ is $58$ in the 319-stem, then we could use $P^{16}$ periodicity to establish that the $h_0$-coexponent of $\frac{\overline{v_1^{64k+96}}}{h_0}$ is 58 for all $k \geq 0$.
\end{rem}

Except for the cases $k = 24$ and $k = 2^j$, we know essentially nothing about the $h_0$-coexponents of $\frac{\overline{v_1^{4k}}}{h_0}$ when $k$ is a multiple of $8$.  For example, when $k = 48$, we cannot even guess the value of the $h_0$-coexponent of $\frac{\overline{v_1^{192}}}{h_0}$ in the $(384-1)$-stem.  This difficulty occurs even though the problem is entirely algebraic in nature, in the sense that it occurs in $H^{**}\cA^{\cl}$, not in stable homotopy.

\subsection{\texorpdfstring{$w_1$}{w1}-periodicity in the cohomology of the \texorpdfstring{$\C$}{C}-motivic Steenrod algebra}
\label{subsctn:w1-algebraic}

In \cref{subsctn:v1-algebraic}, we carefully described some $v_1$-periodic elements in $H^{**}\cA^{\cl}$.  Recall from \cref{Ddi} that $H^{**} \cA^{\cl}$ is isomorphic to the part of $H^{***} \cA$ in Chow degree zero.  Therefore, all of the $v_1$-periodic structure from \cref{subsctn:v1-algebraic} translates to $w_1$-periodic structure in $H^{***} \cA$ in Chow degree zero.

Let $\frac{\overline{w_1^{4k}}}{h_1}$ be the unique non-zero element of $H^{***}\cA$ in degree $(20k-2, 4k, 12k-1)$, which lies in Chow degree zero.  This element corresponds to $\frac{\overline{v_1^{4k}}}{h_0}$ under the isomorphism of \cref{Ddi} from $H^{**} \cA^{\cl}$ to the Chow degree zero part of $H^{***} \cA$.  Consequently, all of the algebraic results in \cref{subsctn:v1-algebraic} can be converted into analogous statements about $\frac{\overline{w_1^{4k}}}{h_1}$.  The notation reflects the fact that $w_1^{4k}$ is the reason that $\frac{\overline{w_1^{4k}}}{h_1}$ is annihilated by $h_1$.

We use the terminology and notation of \cref{subsctn:coexponent} to describe $h_1$-submodules in terms of cogenerators and $h_1$-coexponents, so $\frac{\overline{w_1^{4k}}}{h_1^n}$ is an element $y$ with the property that $h_1^{n-1} \cdot y$ equals $\frac{\overline{w_1^{4k}}}{h_1}$.

\begin{rem}
\label{rem:w1-brackets}
Analogous to \cref{rem:v1-brackets}, the elements $\frac{\overline{w_1^{4k}}}{h_1^n}$ are related to each other via the Massey products that express $w_1$-periodicity.  Specifically, we have
\[
\frac{\overline{w_1^{4j+4k}}}{h_1} = 
\left\langle \frac{\overline{w_1^{4j}}}{h_1}, h_1, \frac{\overline{w_1^{4k}}}{h_1} \right\rangle
\]
with no indeterminacy for all $j$ and $k$.

More generally, assume that $n \leq m$.  For each choice of $\frac{\overline{w_1^{4j}}}{h_1^m}$ and each choice of $\frac{\overline{w_1^{4k}}}{h_1^n}$, the Massey product
\[
\left\langle \frac{\overline{w_1^{4j}}}{h_1^m}, h_1^m, \frac{\overline{w_1^{4k}}}{h_1^n} \right\rangle
\]
has no indeterminacy and equals an element that is a possible choice for $\frac{\overline{w_1^{4j+4k}}}{h_1^n}$.  Beware that the unique element in the Massey product may depend on the choices of $\frac{\overline{w_1^{4j}}}{h_1^m}$ and of $\frac{\overline{w_1^{4k}}}{h_1^n}$.
\end{rem}

For a fixed value of $n$, consider the elements of the form $\frac{\overline{w_1^{4k}}}{h_1^n}$ in $H^{***}\cA$.  Such an element exists only when the $h_1$-coexponent of $\frac{\overline{w_1^{4k}}}{h_1}$ is at least $n$.  This family of existing elements forms a $w_1$-periodic family, but the degree of the periodicity depends on $n$.

\subsection{An additional \texorpdfstring{$v_1$}{v1}-periodic and \texorpdfstring{$w_1$}{w1}-periodic family}
\label{subsctn:v1-algebraic-additional}

We have studied some $v_1$-periodic elements in stems congruent to 7 modulo 8.  There is an additional family of $v_1$-periodic elements in stems congruent to 3 modulo 8.  In a sense, all of these elements fit together into one family.  However, the elements in stems congruent to 3 modulo 8 have slightly different behavior and require slightly different notation to describe.  A similar phenomenon occurs in the classical Adams--Novikov spectral sequence, in which there are differentials and hidden extensions in stems congruent to 3 modulo 8, but not in stems congruent to 7 modulo 8.

Let $\overline{v_1^{4k+2}}$ be the unique non-zero element in $H^{**}\cA^{\cl}$ in degree $(8k+3, 4k+3)$.  These elements are easily identified on a standard Adams chart such as \cite{IWX22b}.  They are the elements at the top of the $(8k+3)$-stem and are usually known as $P^{k} h_0^2 h_2$ or $P^{k} h_1^3$.

The notation reflects the fact that $h_0 v_1^{4k+2}$ is the reason that $\overline{v_1^{4k+2}}$ is annihilated by $h_0$.  Philosophically, $v_1^{4k+2}$ is the reason that $\overline{v_1^{4k+2}}$ is divisible by both $h_0$ and $h_1$.  For example, in the May spectral sequence, $v_1^{4k+2}$ supports a differential whose value is a linear combination of an $h_1$ multiple and an $h_0$ multiple.  In the Cartan--Eilenberg spectral sequence, $v_1^{4k+2}$ supports a differential whose value is an $h_0$ multiple, but this makes room for a hidden $h_0$ extension to an element that is already an $h_1$ multiple.

Each of the elements $\overline{v_1^{4k+2}}$ cogenerates an $h_0$-submodule whose $h_0$-coexponent is $3$.  This $h_0$-submodule is generated by $P^{k} h_2$.

The notation has been arranged carefully for compatibility with $v_1$-periodicity.  This compatibility is most easily seen in \cref{fig:w1-periodic-algebraic}, in which $v_1$-periodic families of elements are arranged in horizontal rows.

Similarly to the $v_1$-periodic elements that were just described, there is an additional family of $w_1$-periodic elements in coweights congruent to 3 modulo 8 that exhibits slightly different behavior, even though all of these elements fit together into one family in a sense.  Let $\overline{w_1^{4k+2}}$ be the unique non-zero element in $H^{***}\cA$ in degree $(20k+9, 4k+3, 12k+6)$.  This element corresponds to $\overline{v_1^{4k+2}}$ under the isomorphism of \cref{Ddi} from $H^{**} \cA^{\cl}$ to the Chow degree zero part of $H^{***} \cA$.  The notation reflects the fact that $h_1 w_1^{4k+2}$ is the reason that $\overline{w_1^{4k+2}}$ is annihilated by $h_1$.  Note that $w_1^{4k+2}$ is the reason that $\overline{w_1^{4k+2}}$ is divisible by both $h_1$ and $h_2$.

\subsection{The structure of \texorpdfstring{$v_1$}{v1}-periodic \texorpdfstring{$h_0$}{h0}-submodules}
\label{subsctn:v1-periodic-h0-submodules}

\cref{table:h0-divisibility} summarizes what we already know about the $v_1$-periodic $h_0$-submodules that are cogenerated by $\frac{\overline{v_1^{4k}}}{h_0}$ and $\overline{v_1^{4k+2}}$.  The table describes the $h_0$-submodules in terms of their cogenerators, $h_0$-coexponents, and generators.  The generator in the last row of the table has no standard name because it lies in a stem within the range of current machine computation but beyond the range of thorough manual exploration.  Translations of the stems and the degrees are displayed because they reveal the numerical patterns more evidently.

\renewcommand{\arraystretch}{2.5} 
\begin{longtable}{lllll}
\caption{Some $v_1$-periodic $h_0$-submodules in the cohomology of the classical Steenrod algebra
\label{table:h0-divisibility}
} \\
\toprule
$\mathrm{stem}+1$ & cogenerator & $\mathrm{degree} + (1,0)$ & $h_0$-coexponent & generator \\
\midrule \endfirsthead
\caption[]{Some $v_1$-periodic $h_0$-submodules in the cohomology of the classical Steenrod algebra}\\
\toprule
$\mathrm{stem}+1$ & cogenerator & $\mathrm{degree} + (1,0)$ & $h_0$-coexponent & generator \\
\midrule \endhead
\bottomrule \endfoot
$2^{j+3}$ & $\displaystyle \frac{\overline{v_1^{2^{j+2}}}}{h_0}$ & $2^{j}\cdot (8,4)$ & $2^{j+2}$ & $h_{j+3}$ \\ 
$8k+4$ & $\overline{v_1^{4k+2}}$ & $k \cdot (8,4) + (4,3)$ & $3$ & $P^{k} h_2$ \\ 
$16k+24$ & $\displaystyle \frac{\overline{v_1^{8k+12}}}{h_0}$ & $(2k+3) \cdot (8,4)$ & $6$ & $P^{2k} i$ \\ 
$32k+48$ & $\displaystyle \frac{\overline{v_1^{16k+24}}}{h_0}$ & $(4k+6) \cdot (8,4)$ & $12$ & $P^{4k} \Delta h_0^2 i$ \\ 
$64k+96$ & $\displaystyle \frac{\overline{v_1^{32k+48}}}{h_0}$ & $(8k + 12) \cdot (8,4)$ & $28$ & $P^{8k} \Delta^3 h_0^2 i$ \\ 
$192$ & $\displaystyle \frac{\overline{v_1^{96}}}{h_0}$ & $24 \cdot (8, 4)$ & $58$ & unnamed \\ 
\end{longtable}
\renewcommand{\arraystretch}{1.0}

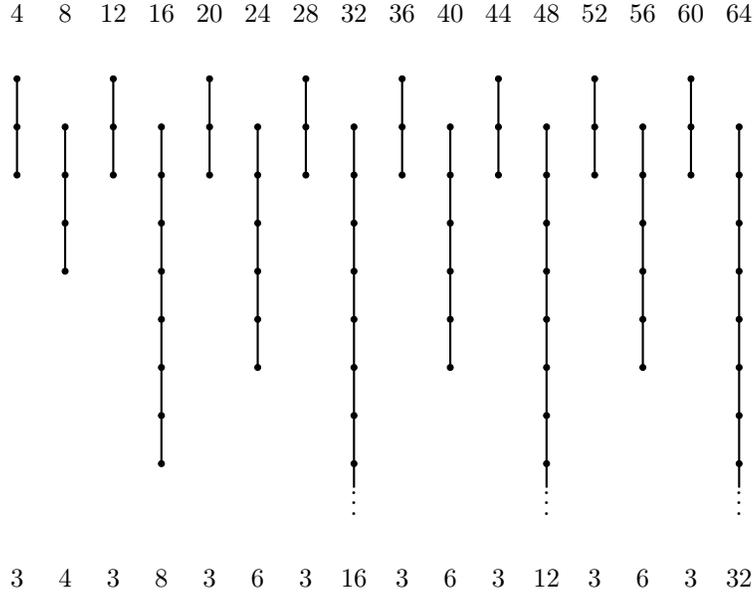
\begin{figure}[H]


\begin{tikzpicture}[scale=0.64]

\pgfmathsetmacro{\minbottom}{100}
\foreach \bottom/\top/\tail/\custombottom in {-1/1/0/0, -3/0/0/0, -1/1/0/0, -7/0/0/0, -1/1/0/0, -5/0/0/0, -1/1/0/0, -3/0/1/16, -1/1/0/0, -5/0/0/0, -1/1/0/0, -3/0/1/12, -1/1/0/0, -5/0/0/0, -1/1/0/0, -3/0/1/32} {
    \ifnum\tail=0
        \ifnum\bottom<\minbottom
            \xdef\minbottom{\bottom} 
        \fi
    \fi
}

\pgfmathsetmacro{\maxtop}{-100}
\foreach \bottom/\top/\tail/\custombottom in {-1/1/0/0, -3/0/0/0, -1/1/0/0, -7/0/0/0, -1/1/0/0, -5/0/0/0, -1/1/0/0, -3/0/1/16, -1/1/0/0, -5/0/0/0, -1/1/0/0, -3/0/1/12, -1/1/0/0, -5/0/0/0, -1/1/0/0, -3/0/1/32} {
    \ifnum\top>\maxtop
        \xdef\maxtop{\top} 
    \fi
}

\foreach \bottom/\top/\tail/\custombottom [count=\i] in {-1/1/0/0, -3/0/0/0, -1/1/0/0, -7/0/0/0, -1/1/0/0, -5/0/0/0, -1/1/0/0, -3/0/1/16, -1/1/0/0, -5/0/0/0, -1/1/0/0, -3/0/1/12, -1/1/0/0, -5/0/0/0, -1/1/0/0, -3/0/1/32} {

    \pgfmathsetmacro{\topLabel}{int(4*\i)} 
    \node[above] at (\i,  \maxtop + 1) {\textnormal{\topLabel}}; 
    
    \ifnum\tail=1
        \node[below] at (\i, \minbottom - 2) {\custombottom}; 
    \else
        \pgfmathsetmacro{\bottomLabel}{int(\top-\bottom + 1)} 
        \node[below] at (\i, \minbottom - 2) {\bottomLabel}; 
    \fi

    \draw[thick] (\i, \bottom) -- (\i, \top);
    
    \pgfmathsetmacro{\n}{\top-\bottom} 
    
    \foreach \k in {0,...,\n} {
        \pgfmathsetmacro{\y}{\bottom + \k}
        \filldraw (\i, \y) circle (1.75pt); 
    }
    
    \ifnum\tail=1
        \draw[thick] (\i, \minbottom) -- (\i, \bottom);
        \pgfmathsetmacro{\m}{\bottom-\minbottom -1}
        \foreach \k in {0,...,\m} {
            \pgfmathsetmacro{\y}{\minbottom + \k}
            \filldraw (\i, \y) circle (1.75pt); 
        }
        \draw[thick] (\i, \minbottom) -- (\i, \minbottom - 0.5); 
        
        \node at (\i, \minbottom - 0.65) {\scalebox{1}{\vdots}};
    \fi

}



\end{tikzpicture}

\caption{Some $v_1$-periodic $h_0$-submodules in the classical Adams $E_2$-page, or some $w_1$-periodic $h_1$-submodules in the $\C$-motivic Adams $E_2$-page.  Numbers along the top are $(\mathrm{stem}+ 1)$, resp., $(\mathrm{coweight} + 1)$. Numbers along the bottom are coexponents.}
\label{fig:w1-periodic-algebraic}
\end{figure}

Some of the information in \cref{table:h0-divisibility} is displayed graphically in \cref{fig:w1-periodic-algebraic}.  The figure shows the cyclic $h_0$-submodules of the classical Adams $E_2$-page that are cogenerated by $\frac{\overline{v_1^{4k}}}{h_0}$ and $\overline{v_1^{4k+2}}$. The numbers along the top of the figure indicate the stem plus one. The numbers along the bottom of the figure indicate $h_0$-coexponents.  The elements are arranged such that every horizontal row is a $v_1$-periodic family. The element $\frac{\overline{v_1^{4k}}}{h_0}$ lies at the top of column $8k$, and the element $\overline{v_1^{4k+2}}$ lies at the top of column $8k+4$. Beware that the vertical alignment of \cref{fig:w1-periodic-algebraic} does not reflect the Adams filtration.  On a standard Adams chart, $v_1$-periodic families lie along lines of slope $\frac{1}{2}$, not along horizontal lines. The vertical alignments of the columns indicate why the stems congruent to 7 modulo 8 must be studied somewhat separately from the stems congruent to 3 modulo 8.  

For an explanation of the secondary caption on \cref{fig:w1-periodic-algebraic}, see the end of \cref{subsctn:w1-periodic-h1-submodules} below.

The data in \cref{table:h0-divisibility} suggests the following conjecture.

\begin{conj}
\label{conj:h0-size}
For each fixed value of $j \geq 3$, the $h_0$-submodules cogenerated by $\displaystyle \frac{\overline{v_1^{3 \cdot 2^j + i \cdot 2^{j+1}}}}{h_0}$ have the same $h_0$-coexponent for all $i$.
\end{conj}

In number-theoretic language, \cref{conj:h0-size} says that the  $h_0$-coexponent of the $h_0$-submodule cogenerated by $\frac{\overline{v_1^{4k}}}{h_0}$ depends only on the  2-adic valuation of $4k$ (i.e., on how many powers of $2$ divide $4k$), except when $4k$ is a power of $2$.

The conjecture is already verified for $3 \leq j \leq 5$ in \cref{lem:h0-coexponents}.  The difficulty starting with $j = 6$ is discussed in \cref{rem:h0-coexponent-96}.

Assuming \cref{conj:h0-size}, we can then reduce the problem to a more specific question.

\begin{Que}
\label{question:h0-size}
For $j \geq 3$, what is the $h_0$-coexponent of $\frac{\overline{v_1^{3 \cdot 2^j}}}{h_0}$?
\end{Que}

The first answers to \cref{question:h0-size} are 6, 12, 28, and 58. We have little idea about the later values, except to guess that each value is approximately twice the previous value.

\subsection{The structure of \texorpdfstring{$w_1$}{w1}-periodic \texorpdfstring{$h_1$}{h1}-submodules}
\label{subsctn:w1-periodic-h1-submodules}

The isomorphism of \cref{Ddi} from $H^{**} \cA^{\cl}$ to the part of $H^{***} \cA$ in Chow degree zero transfers all of the structure about $h_0$-submodules in \cref{subsctn:v1-periodic-h0-submodules} to structure about $h_1$-submodules.

\begin{prop}
\cref{table:h1-divisibility} describes some cyclic $h_1$-submodules of $H^{***}\cA$ in Chow degree zero.
\end{prop}

\begin{proof}
We use the isomorphism of \cref{Ddi} from $H^{**} \cA^{\cl}$ to the Chow degree zero part of $H^{***} \cA$.  Apply this isomorphism to \cref{lem:h0-coexponent-2^j}, \cref{lem:h0-coexponents}, and \cref{lem:h0-coexponent-96}, or directly to \cref{table:h0-divisibility}.
\end{proof}

\renewcommand{\arraystretch}{2.0}
\begin{longtable}{lllll}
\caption{Some $w_1$-periodic $h_1$-submodules in the cohomology of the $\C$-motivic Steenrod algebra
\label{table:h1-divisibility}
} \\
\toprule
$\mathrm{coweight} + 1$ & $\mathrm{degree} + (2,0,1)$ & cogenerator  & $h_1$-coexponent & generator \\ 
\midrule \endfirsthead
\caption[]{Some $w_1$-periodic $h_1$-submodules in the cohomology of the $\C$-motivic Steenrod algebra}\\
\toprule
$\mathrm{coweight} + 1$ & $\mathrm{degree} + (2,0,1)$ & cogenerator & $h_1$-coexp. & generator \\ 
\midrule \endhead
\bottomrule \endfoot
$2^{j+3}$ & $2^{j} \cdot (20, 4, 12)$ & $\frac{\overline{w_1^{2^{j+2}}}}{h_1}$ & $2^{j+2}$ & $h_{j+4}$ \\ 
$8k+4$ & $k \cdot (20, 4, 12) + (11, 3, 7)$ & $\overline{w_1^{4k+2}}$ & $3$ & $h_3 g^{k}$ \\ 
$16k+24$ & $(2k + 3) \cdot (20, 4, 12)$ & $\frac{\overline{w_1^{8k+12}}}{h_1}$ & $6$ & $g^{2k} i_1$ \\
$32k+48$ & $(4k + 6) \cdot (20, 4, 12)$ & $\frac{\overline{w_1^{16k+24}}}{h_1}$ & $12$ & $g^{4k} \Delta_1 h_1^2 i_1$ \\
$64k+96$ & $(8k + 12) \cdot (20, 4, 12)$ & $\frac{\overline{w_1^{32k+48}}}{h_1}$ & $28$ & $g^{8k} \Delta_1^3 h_1^2 i_1$ \\
$192$ & $24 \cdot (20, 4, 12)$ & $\frac{\overline{w_1^{96}}}{h_1}$ & $58$ & unnamed \\ 
\end{longtable}
\renewcommand{\arraystretch}{1.0}

\cref{table:h1-divisibility} summarizes what we already know about the $w_1$-periodic $h_1$-submodules that are cogenerated by  $\frac{\overline{w_1^{4k}}}{h_1}$ and $\overline{w_1^{4k+2}}$. The table describes the $h_1$-submodules in terms of their cogenerators, $h_1$-coexponents, and generators. Here we use the coweight, which equals the stem minus the motivic weight. Note that $0$ is the coweight of $h_1$, so every cyclic $h_1$-submodule  is concentrated in a single coweight.  This explains why coweight is useful for describing the degree of a cyclic $h_1$-submodule. Translations of the coweights and the degrees are displayed because they reveal the numerical patterns more evidently.  Beware that the table gives the degrees of the cogenerators, but not of the generators.

Except for the cases $k = 24$ and $k = 2^j$, we know essentially nothing about the $h_1$-coexponents of $\frac{\overline{w_1^{4k}}}{h_1}$ when $k$ is a multiple of $8$.  For example, when $k = 48$, we cannot even guess the value of the $h_1$-coexponent of $\frac{\overline{w_1^{192}}}{h_1}$ in coweight $(384-1)$.  This difficulty occurs in precisely the same way that it does for $v_1$-periodic $h_0$-coexponents at the end of \cref{subsctn:v1-algebraic}.

Under suitable interpretations, \cref{fig:w1-periodic-algebraic} displays both the $v_1$-periodic $h_0$-submodules of $H^{**} \cA^{\cl}$ and also the $w_1$-periodic $h_1$-submodules of $H^{***} \cA$.  In the latter interpretation,  the $w_1$-periodic families of elements are arranged in horizontal rows. Moreover, the vertical towers represent $h_1$-submodules, which appear along lines of slope $1$ in a standard $\C$-motivic Adams chart.

\cref{conj:h0-size} has an analogue for $h_1$-submodules.  Namely, we guess that for each fixed value of $j \geq 3$, the $h_1$-submodules cogenerated by $\displaystyle \frac{\overline{w_1^{3 \cdot 2^j + i \cdot 2^{j+1}}}}{h_1}$ have the same $h_1$-coexponent for all $i$. The data in \cref{table:h1-divisibility} supports this prediction. In number-theoretic language, we expect that the $h_1$-coexponent of the $h_1$-submodule cogenerated by $\frac{\overline{w_1^{4k}}}{h_1}$ depends only on the 2-adic valuation of $4k$ (i.e., on how many powers of $2$ divide $4k$), except when $4k$ is a power of $2$. More specifically, for $j \geq 3$, we ask for the $h_1$-coexponent of the $h_1$-submodule cogenerated by $\frac{\overline{w_1^{3 \cdot 2^j}}}{h_1}$.

\section{Classical \texorpdfstring{$v_1$}{v1}-periodicity and \texorpdfstring{$\C$}{C}-motivic \texorpdfstring{$w_1$}{w1}-periodicity in homotopy}
\label{sctn:v1-w1-periodic-homotopy}

\subsection{Classical \texorpdfstring{$v_1$}{v1}-periodic homotopy}
\label{subsctn:v1-homotopy}
We recall the number-theoretically interesting behavior of the $2$-primary $v_1$-periodic stable homotopy groups in stems congruent to 3 and 7 modulo 8.  There is non-zero $2$-primary $v_1$-periodic stable homotopy in stems congruent to 0, 1, and 2 modulo 8, but we exclude those groups from our discussion because they are less interesting from a number-theoretic perspective.

Starting with the 3-stem, the $2$-primary $v_1$-periodic stable homotopy groups are cyclic and have exponents
\begin{equation}
\label{eq:v1-exponent}
3, 4, 3, 5, 3, 4, 3, 6, 3, 4, 3, 5, 3, 4, 3, 7, \ldots
\end{equation}
Here, the exponent is the smallest number $e$ such that $2^e$ annihilates the group.  This sequence displays a regular number-theoretic pattern.  Namely,
the $k$th entry is equal to the $2$-adic valuation of the integer $8k$.  In other words, the order of the $(4k-1)$-stem depends only on the $2$-adic valuation of $k$.  In even more concrete terms, every other group is of exponent $3$, every fourth group is of exponent $4$, every eighth group is of
exponent $5$, and so forth.  \cref{fig:v1-periodic} illustrates the structure.

We briefly summarize the history of these $v_1$-periodic computations.  Milnor and Kervaire \cite{MK60} proved that the values in sequence \eqref{eq:v1-exponent} are lower bounds.  Adams \cite{Ada66a} showed that the exponent is at most one more than the values in the sequence.  Quillen \cite[Theorem 1.1]{Qui71} and Sullivan \cite{Sul74} independently proved the Adams Conjecture, which confirmed that the exponents are as shown in sequence \eqref{eq:v1-exponent}.

\subsection{Classical \texorpdfstring{$v_1$}{v1}-periodic homotopy in the Adams \texorpdfstring{$E_\infty$}{E infinity}-page}
\label{subsctn:v1-periodic-homotopy-Adams}

The elements $\frac{\overline{v_1^{4k}}}{h_0}$ and $\overline{v_1^{4k+2}}$ are all non-zero permanent cycles in the classical Adams $E_\infty$-page.  They detect the unique $2$-torsion $v_1$-periodic elements in stems congruent to 7 modulo 8 and 3 modulo 8 respectively.  The $h_0$-coexponents of $\frac{\overline{v_1^{4k}}}{h_0}$ in the $E_2$-page are larger than their $h_0$-coexponents in the $E_\infty$-page because there are Adams differentials on many of the elements that are cogenerated by $\frac{\overline{v_1^{4k}}}{h_0}$.

In both the $E_2$-page and the $E_\infty$-page, the $h_0$-coexponents of $\frac{\overline{v_1^{4k}}}{h_0}$ display number-theoretic structure in the sense that they depend on the 2-adic valuations of $4k$.  However, the \emph{values} of these $h_0$-coexponents are very different in the $E_2$-page than in the $E_\infty$-page.  

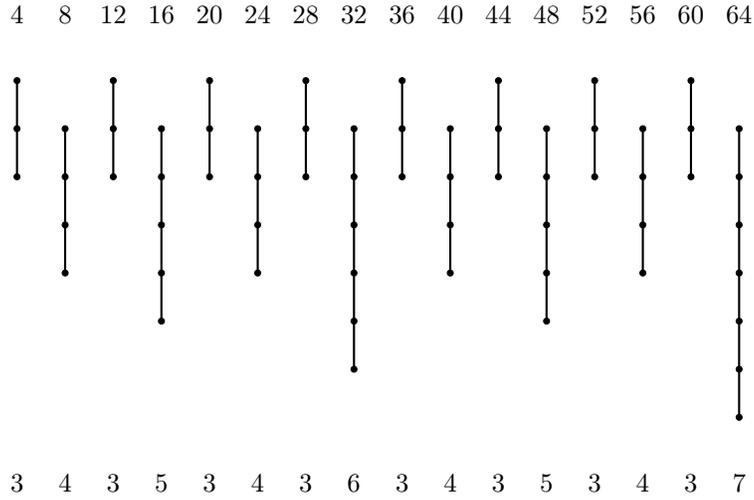
\begin{figure}[H]


\begin{tikzpicture}[scale=0.64]

\pgfmathsetmacro{\minbottom}{100}
\foreach \bottom/\top/\tail/\custombottom in {-1/1/0/0, -3/0/0/0, -1/1/0/0, -4/0/0/0, -1/1/0/0, -3/0/0/0, -1/1/0/0, -5/0/0/0, -1/1/0/0, -3/0/0/0, -1/1/0/0, -4/0/0/0, -1/1/0/0, -3/0/0/0, -1/1/0/0, -6/0/0/0} {
    \ifnum\tail=0
        \ifnum\bottom<\minbottom
            \xdef\minbottom{\bottom} 
        \fi
    \fi
}

\pgfmathsetmacro{\maxtop}{-100}
\foreach \bottom/\top/\tail/\custombottom in {-1/1/0/0, -3/0/0/0, -1/1/0/0, -4/0/0/0, -1/1/0/0, -3/0/0/0, -1/1/0/0, -5/0/0/0, -1/1/0/0, -3/0/0/0, -1/1/0/0, -4/0/0/0, -1/1/0/0, -3/0/0/0, -1/1/0/0, -6/0/0/0} {
    \ifnum\top>\maxtop
        \xdef\maxtop{\top} 
    \fi
}

\foreach \bottom/\top/\tail/\custombottom [count=\i] in {-1/1/0/0, -3/0/0/0, -1/1/0/0, -4/0/0/0, -1/1/0/0, -3/0/0/0, -1/1/0/0, -5/0/0/0, -1/1/0/0, -3/0/0/0, -1/1/0/0, -4/0/0/0, -1/1/0/0, -3/0/0/0, -1/1/0/0, -6/0/0/0} {

    \pgfmathsetmacro{\topLabel}{int(4*\i)} 
    \node[above] at (\i,  \maxtop + 1) {\textnormal{\topLabel}}; 

    \ifnum\tail=1
        \node[below] at (\i, -2) {\custombottom}; 
    \else
        \pgfmathsetmacro{\bottomLabel}{int(\top-\bottom + 1)} 
        \node[below] at (\i, \minbottom - 1) {\bottomLabel}; 
    \fi

    \draw[thick] (\i, \bottom) -- (\i, \top);
    
    \pgfmathsetmacro{\n}{\top-\bottom} 
    
    \foreach \k in {0,...,\n} {
        \pgfmathsetmacro{\y}{\bottom + \k}
        \filldraw (\i, \y) circle (1.75pt); 
    }
    
    \ifnum\tail=1
        \draw[thick] (\i, \minbottom) -- (\i, \bottom);
        \pgfmathsetmacro{\m}{\bottom-\minbottom - 1}
        \foreach \k in {0,...,\m} {
            \pgfmathsetmacro{\y}{\minbottom + \k}
            \filldraw (\i, \y) circle (1.75pt); 
        }
        \draw[thick] (\i, \minbottom) -- (\i, \minbottom - 0.5); 
        
        \node at (\i, \minbottom - 0.5) {\scalebox{1.5}{\vdots}};
    \fi
    }
\end{tikzpicture}

\caption{Some $v_1$-periodic $h_0$-submodules in the classical Adams $E_\infty$-page. Numbers along the top are $(\mathrm{stem}+ 1)$. Numbers along the bottom are coexponents.}
\label{fig:v1-periodic}
\end{figure}

\cref{fig:v1-periodic} displays the cyclic $h_0$-submodules of the classical Adams $E_\infty$-page that are cogenerated by $\frac{\overline{v_1^{4k}}}{h_0}$ and $\overline{v_1^{4k+2}}$.  The figure is qualitatively similar to \cref{fig:w1-periodic-algebraic}, but the $h_0$-coexponents of the $h_0$-submodules take different values.

The numbers along the top of the figure indicate the stem plus one.  The numbers along the bottom of the figure indicate $h_0$-coexponents.  As in \cref{fig:w1-periodic-algebraic}, the elements are arranged such that every
horizontal row is a $v_1$-periodic family.  The element $\frac{\overline{v_1^{4k}}}{h_0}$ lies at the top of column $8k$, and the element $\overline{v_1^{4k+2}}$ lies at the top of column $8k+4$.  Beware that the vertical alignment of \cref{fig:v1-periodic} does not reflect the Adams filtration.  On a standard Adams chart, $v_1$-periodic families lie along lines of slope $\frac{1}{2}$, not along horizontal lines.  As in \cref{fig:w1-periodic-algebraic}, the vertical alignments of the columns indicate why the stems congruent to 7 modulo 8 must be studied somewhat separately from the stems congruent to 3 modulo 8.

\subsection{\texorpdfstring{$w_1$}{w1}-periodic homotopy}
\label{subsctn:w1-periodic-homotopy}

In \cref{subsctn:w1-periodic-h1-submodules}, we studied some $h_1$-submodules in the $\C$-motivic Adams $E_2$-page.  We now push further from algebraic phenomena into results about $\C$-motivic homotopy.

There are several complications that can occur in the transition from algebra to homotopy.  There turn out to be no problems with hidden extensions, but we do need to worry about the presence of Adams differentials.  Our results in \cref{sctn:Adams-differentials} give some information about Adams differentials, but we do not know the values of all Adams differentials on the $w_1$-periodic $E_2$-page elements that are described in \cref{subsctn:w1-periodic-h1-submodules}.

\begin{lem}
\label{prop:no-differentials}
No element of the form $\frac{\overline{w_1^{4k}}}{h_1^n}$ or $\frac{\overline{w_1^{4k+2}}}{h_1^n}$ is hit by an Adams differential.
\end{lem}

\begin{proof}
All of these elements have Chow degree zero.  The Adams $d_r$ differential increases the Chow degree by $r-1$.  The $\C$-motivic Adams $E_2$-page is zero in negative Chow degrees by \cref{prop:Chow-degree-nonneg}, so all differentials entering Chow degree zero are null.
\end{proof}

\begin{lem}
\label{prop:no-eta-extns-in}
No element of the form $\frac{\overline{w_1^{4k}}}{h_1^n}$ or $\frac{\overline{w_1^{4k+2}}}{h_1^n}$ is the target of a hidden $\eta$-extension.
\end{lem}

\begin{proof}
Each element under consideration has Chow degree zero, so it has a degree of the form $\left( s,f,\frac{s+f}{2} \right)$ for some $s$ and $f$.  If that element were the target of a hidden $\eta$ extension, then the source of the extension would have degree $\left( s-1, f', \frac{s+f}{2} -1 \right)$, where $f' < f-1$.  Then the Chow degree of this source would be $f' - f  + 1$,  which is negative.  The $\C$-motivic Adams $E_2$-page is zero in negative Chow degrees by \cref{prop:Chow-degree-nonneg}.
\end{proof}

\begin{lem}
\label{prop:no-tau-extns-in}
No element of the form $\frac{\overline{w_1^{4k}}}{h_1^n}$ or $\frac{\overline{w_1^{4k+2}}}{h_1^n}$ is the target of a (hidden or non-hidden) $\tau$-extension.
\end{lem}

\begin{proof}
Each element under consideration has Chow degree zero, so it has a degree of the form $\left( s,f,\frac{s+f}{2} \right)$ for some $s$ and $f$.  If that element were the target of a $\tau$ extension, then the source of the extension would have degree $\left( s, f', \frac{s+f}{2} +1 \right)$, where $f' \leq f$.  (Note that $f' = f$ precisely when the $\tau$-extension is not hidden.)  Then the Chow degree of this source would be $f' - f - 2$, which is negative.  The $\C$-motivic Adams $E_2$-page is zero in negative Chow degrees by \cref{prop:Chow-degree-nonneg}.
\end{proof}

\begin{prop}
\label{prop:h1hn}
The element $h_1 h_n$ is a permanent cycle in the $\C$-motivic Adams spectral sequence.
\end{prop}

\begin{proof}
Using \cite{Mahowald77}, we know that $h_1 h_n$ is a permanent cycle in the $\tau$-localization of the $\C$-motivic Adams spectral sequence.  It is possible that $\tau^m h_1 h_n$ supports differentials for small values of $m$, but there exists a smallest $N$ such that $\tau^N h_1 h_n$ is a permanent cycle.  Let $\alpha$ be the element in $\pi_{2^n, 2^{n-1}+1-N}$ that is detected by $\tau^N h_1 h_n$.

The framework of \cite[Chapter 6]{Isa19} tells us that if an element of $\pi_{s,w}$ is not divisible by $\tau$, then it is detected in filtration $2w-s$ in the Adams--Novikov spectral sequence.  In particular, $\alpha$ is detected in filtration $2-2N$ in the Adams--Novikov spectral sequence.  The only possibility is that $N = 0$.
\end{proof}

\begin{rem}
\cref{prop:h1hn} is not an independent proof of the existence of Mahowald's $\eta_n$ family.  Rather, we use Mahowald's classical result to extend to the $\C$-motivic context.
\end{rem}

We use the terminology and notation of \cref{subsctn:coexponent} to describe $\eta$-submodules of $\pi_{*,*}$ in terms of cogenerators and $\eta$-coexponents.  Beware that the generator of an $\eta$-submodule is not necessarily annihilated by $2$.  In the language of \cref{subsctn:coexponent}, the characteristic of an $\eta$-submodule can be larger than $2$.  For example, the element $\eta_6$ in $\pi_{64,33}$ is detected by $h_1 h_6$, and it supports a $2$ extension that is hidden in the $\C$-motivic Adams $E_\infty$-page \cite{IWX20}.  The element $\eta_6$ is a generator of one of the explicit $\eta$-submodules that we study.    We cannot control these hidden $2$ extensions on the generators of our $\eta$-submodules, so we do not discuss characteristics of $\eta$-submodules in detail.

\begin{thm}
\label{thm:C-homotopy}
For all $n \geq 1$, there are non-zero elements $\gamma_n$ in $\pi_{20n-2,12n-1}$ such that:
\begin{enumerate}
\item 
$\gamma_n$ is detected by $\frac{\overline{w_1^{4n}}}{h_1}$ in filtration $4n$ in the $\C$-motivic Adams spectral sequence.
\item
$\gamma_n$ is annihilated by $\eta$.
\item
the elements form a $w_1$-periodic family in the sense that $\gamma_{n+1}$ is contained in the Toda bracket $\langle \gamma_n, \eta, \gamma_1 \rangle$ for all $n \geq 1$.
\item
$\gamma_{2^j}$ equals $\eta^{2^{j+2}-2} \cdot \eta_{j+4}$, where $\eta_{j+4}$ is the Mahowald element in the $2^{j+4}$-stem.
\item 
the $\eta$-coexponent of $\gamma_n$ is given in \cref{table:eta-size} for some values of $n$.
\end{enumerate}
\end{thm}

\begin{proof}
The proofs of (1), (2), and (3) proceed by induction.  Let $\gamma_1$ be an element of $\pi_{18,11}$ that is detected by $\frac{\overline{w_1^4}}{h_1}$. We already know that $\gamma_1$ is annihilated by $\eta$ \cite{IWX20}.

Inductively define $\gamma_{n+1}$ to be an element in the Toda bracket
\[
\left\langle
\gamma_n, \eta, \gamma_1
\right\rangle.
\]
The bracket is well-defined by the induction hypothesis.  The Moss convergence theorem \cite{Mos70} then tells us that $\gamma_{n+1}$ is detected by
\[
\frac{\overline{w_1^{4n+4}}}{h_1} = 
\left\langle
\frac{\overline{w_1^{4n}}}{h_1}, h_1, \frac{\overline{w_1^4}}{h_1} 
\right\rangle.
\]
Beware that the Moss convergence theorem requires that there are no crossing Adams differentials.  In this case, there are no crossing differentials because the sources of such differentials would have negative Chow degree, and there are no non-zero elements in negative Chow degree by \cref{prop:Chow-degree-nonneg}.

To show that $\gamma_{n+1}$ is annihilated by $\eta$, consider the shuffle
\[
\eta \langle \gamma_n, \eta, \gamma_1 \rangle =
\langle \eta, \gamma_n, \eta \rangle \gamma_1.
\]
This last expression equals $\nu \gamma_n \gamma_1$ by \cite[Theorem 3.6]{Toda62}, which equals zero because $\nu \gamma_1$ is zero \cite{IWX20}.

For part (4), we use \cref{prop:h1hn} together with the observation that $\frac{\overline{w_1^{2^{j+2}}}}{h_1}$ equals $h_1^{2^{j+2}-2} \cdot h_1 h_{j+4}$.  This last claim is a translation along \cref{Ddi} of the observation that $\frac{\overline{v_1^{2^{j+2}}}}{h_0}$ equals $h_0^{2^{j+2}-2} \cdot h_0 h_{j+3}$ in the classical Adams $E_2$-page.

Finally, we must study the $\eta$-coexponent of $\gamma_n$.  \cref{prop:h1hn} shows that $h_1 h_{j+1}$ is a permanent cycle.  Therefore, $\frac{\overline{w_1^{2^{j-1}}}}{h_1}$ has $h_1$-coexponent equal to $2^{j-1} - 1$ in the Adams $E_\infty$-page.  \cref{prop:no-eta-extns-in} rules out hidden $\eta$ extensions that could increase this coexponent in homotopy.

Similarly, \cref{i1g^2k} and \cref{Delta1h1^2i1g^4k}, to be proved later, show that $\frac{\overline{w_1^{8k+12}}}{h_1}$ and $\frac{\overline{w_1^{16k+24}}}{h_1}$ have $h_1$-coexponents equal to $6$ and $12$ respectively in the Adams $E_\infty$-page.  \cref{prop:no-eta-extns-in} rules out hidden $\eta$ extensions that could increase this coexponent in homotopy.
\end{proof}

\begin{rem}
\label{rem:C-homomtopy-noncircular}
The proofs of parts (1), (2), (3), and (4) of \cref{thm:C-homotopy} are self-contained.  Only the proof of part (5) depends on results that appear later in this manuscript.  We emphasize this point to clarify the absence of circular arguments.  Also, parts (1), (2), and (3) were previously proved by Andrews \cite{And}.
\end{rem}

\begin{rem}
Beware that the Toda brackets in part (3) of \cref{thm:C-homotopy} may have indeterminacies.  
\end{rem}

\begin{rem}
\label{rem:unknown-coexponents}
In part (5) of \cref{thm:C-homotopy}, the $\eta$-coexponent of $\gamma_n$ is unknown for many values of $n$.  The smallest such value is $n = 12$.  There are two potential sources for this uncertainty.  First, we have incomplete knowledge about $h_1$-coexponents in the $\C$-motivic Adams $E_2$-page, as discussed in \cref{subsctn:w1-periodic-h1-submodules}.  Second, even if we did know more about the Adams $E_2$-page, it is possible that the presence of non-zero of Adams differentials could affect the $\eta$-coexponent of $\gamma_n$.
\end{rem}

\cref{table:eta-size} describes some $\eta$-submodules of $\pi_{*,*}$ in terms of their cogenerators, $\eta$-coexponents, and generators.  Translations of the stems and the degrees are displayed because they reveal the numerical patterns more evidently.  The square brackets in the last column refer to elements in homotopy that are detected by the Adams $E_\infty$-page elements that are shown inside the brackets.  These specific detecting elements arise explicitly in \cref{sctn:Adams-differentials}.

\begin{longtable}{llllll}
\caption{Some $\eta$-submodules in $\C$-motivic stable homotopy
\label{table:eta-size}}
\\
\toprule
coweight & cogenerator & $\mathrm{degree} + (2,1)$ & $\eta$-coexp- & generator & proof \\
$+ 1$ & & & onent & & \\
\midrule \endfirsthead
\caption[]{Some $\eta$-submodules in $\C$-motivic stable homotopy}\\
\toprule
coweight & cogenerator & $\mathrm{degree} + (2,1)$ & $\eta$-coexp- & generator & proof \\
$+ 1$ & & & onent & & \\

\midrule \endhead
\bottomrule \endfoot
$2^{j+3}$ & $\gamma_{2^j}$ & $2^{j} \cdot (20, 12)$ & $2^{j+2} - 1$ & $\eta_{j+4}$ & Prop.\ \ref{prop:h1hn} \\
$8n+4$ & $\delta_n$ ($n\geq 1$) & $n \cdot (20,12) + (11, 7)$ & $2$ & $[h_1 h_3 g^n]$ & Prop.\ \ref{h1h3g^k} \\
$16k+24$ & $\gamma_{2k+3}$ & $(2k+3) \cdot (20, 12)$ & $6$ & $[g^{2k} i_1]$ & Prop.\ \ref{i1g^2k}  \\
$32k+48$ & $\gamma_{4k+6}$ & $(4k+6) \cdot (20, 12)$ & $12$ & $[g^{4k} \Delta_1 h_1^2 i_1]$ & Prop.\ \ref{Delta1h1^2i1g^4k} \\
\end{longtable}

\begin{thm}
\label{thm:C-homotopy-beta}
For all $n \geq 1$, there are non-zero elements $\delta_n$ in $\pi_{20n+9,12n+6}$ such that:
\begin{enumerate}
\item 
$\delta_n$ is detected by $\overline{w_1^{4n+2}}$ in filtration $4n+3$ in the $\C$-motivic Adams spectral sequence.
\item
$\delta_n$ is annihilated by $\eta$.
\item
the elements form a $w_1$-periodic family in the sense that $\delta_{n+1}$ is contained in the Toda bracket $\langle \delta_n, \eta, \gamma_1 \rangle$ for all $n$.
\item 
the $\eta$-coexponent of $\delta_n$ is 2.
\end{enumerate}
\end{thm}

\begin{proof}
Identically to the proof of \cref{thm:C-homotopy}, the proofs of (1), (2), and (3) proceed by induction.  Let $\delta_1$ be an element of $\pi_{29,18}$ that is detected by $\overline{w_1^6}$.  We already know that $\delta_1$ is annihilated by $\eta$ \cite{IWX20}.

Inductively define $\delta_{n+1}$ to be an element in the Toda bracket
\[
\left\langle
\delta_n, \eta, \gamma_1 \right\rangle.
\]
The bracket is well-defined by the inductive hypothesis.  The Moss convergence theorem \cite{Mos70} then tells us that $\delta_{n+1}$ is detected by
\[
\overline{w_1^{4n+6}} = 
\left\langle
\overline{w_1^{4n+2}}, h_1, \frac{\overline{w_1^4}}{h_1} \right\rangle.
\]
Beware that the Moss convergence theorem requires that there are no crossing Adams differentials.  In this case, there are no crossing differentials because the source of such differentials would have a negative Chow degree.

To show that $\delta_{n+1}$ is annihilated by $\eta$, consider the shuffle
\[
\eta \langle \delta_n, \eta, \gamma_1 \rangle =
\langle \eta, \delta_n, \eta \rangle \gamma_1.
\]
This last expression equals $\nu \delta_n \gamma_1$ by \cite[Theorem 3.6]{Toda62}, which equals zero because $\nu \gamma_1$ is zero \cite{IWX20}.

Finally, we must study the $\eta$-coexponent of $\delta_n$.  \cref{Cd2} and \cref{h1h3g^k} show that $h_1 h_3 g^k$ is a permanent cycle, but $h_3 g^k$ is not.  Therefore, $\overline{w_1^{4n+2}}$ has $h_1$-coexponent equal to $2$ in the Adams $E_\infty$-page.  \cref{prop:no-eta-extns-in} rules out hidden $\eta$ extensions that could increase this coexponent in homotopy.
\end{proof}

\begin{rem}
\label{rem:C-homomtopy-beta-noncircular}
The proofs of parts (1), (2), and (3) of \cref{thm:C-homotopy-beta} are self-contained.  Only the proof of part (4) depends on results that appear later in this manuscript.  We emphasize this point to clarify the absence of circular arguments.  Also, parts (1), (2), and (3) were previously proved by Andrews \cite{And}.
\end{rem}

\begin{rem}
In \cref{thm:C-homotopy-beta}, one can define $\delta_0$ to be an element in $\pi_{9,6}$ that is detected by $\overline{w_1^2}$, which is more commonly known as $h_1^2 h_3$.  However, this element is qualitatively different.  Because of the presence of $h_1 c_0$ in higher filtration, there is a hidden $\eta$ extension, and $\delta_0$ is not annihilated by any power of $\eta$.  Rather, in coweight $3$, the element $\sigma$ of $\pi_{7,4}$ generates the $\eta$-submodule $\Z[\eta]/(16, 2 \eta)$, and $\delta_0$ equals $\eta^2 \sigma$.  In a sense, the $\eta$-coexponent of $\delta_0$ equals $3$ because it can be divided twice by $\eta$.  However, the product $\eta \delta_0$ is not zero, so $\delta_0$ is not a cogenerator.

Also, $\delta_1$ cannot be constructed as the Toda bracket $\left\langle \delta_0, \eta, \gamma_1 \right\rangle$ because the bracket is not well-defined.
\end{rem}

\begin{figure}[H]


\begin{tikzpicture}[scale=0.64]

\pgfmathsetmacro{\minbottom}{100}
\foreach \bottom/\top/\tip/\tail/\custombottom in {-1/1/1/0/inf, -2/0/0/0/0, 0/1/0/0/0, -6/0/0/0/0, 0/1/0/0/0, -5/0/0/0/0, 0/1/0/0/0, -3/0/0/1/15, 0/1/0/0/0, -5/0/0/0/0, 0/1/0/0/0, -3/0/0/1/12, 0/1/0/0/0, -5/0/0/0/0, 0/1/0/0/0, -3/0/0/1/31} {
    \ifnum\tail=0
        \ifnum\bottom<\minbottom
            \xdef\minbottom{\bottom} 
        \fi
    \fi
}

\pgfmathsetmacro{\maxtop}{-100}
\foreach \bottom/\top/\tip/\tail/\custombottom in {-1/1/1/0/inf, -2/0/0/0/0, 0/1/0/0/0, -6/0/0/0/0, 0/1/0/0/0, -5/0/0/0/0, 0/1/0/0/0, -3/0/0/1/15, 0/1/0/0/0, -5/0/0/0/0, 0/1/0/0/0, -3/0/0/1/12, 0/1/0/0/0, -5/0/0/0/0, 0/1/0/0/0, -3/0/0/1/31} {
    \ifnum\top>\maxtop
        \xdef\maxtop{\top} 
    \fi
}

\foreach \bottom/\top/\tip/\tail/\custombottom [count=\i] in {-1/1/1/0/inf, -2/0/0/0/0, 0/1/0/0/0, -6/0/0/0/0, 0/1/0/0/0, -5/0/0/0/0, 0/1/0/0/0, -3/0/0/1/15, 0/1/0/0/0, -5/0/0/0/0, 0/1/0/0/0, -3/0/0/1/12, 0/1/0/0/0, -5/0/0/0/0, 0/1/0/0/0, -3/0/0/1/31} {

    \pgfmathsetmacro{\topLabel}{int(4*\i)} 
    \node[above] at (\i,  \maxtop + 1) {\textnormal{\topLabel}}; 
    
    \ifnum\tail=1
        \node[below] at (\i, \minbottom-2) {\custombottom};  
    \else
        \ifnum0=\pdfstrcmp{\custombottom}{inf} 
            \node[below] at (\i, \minbottom-2) {$\infty$};  
        \else
            \pgfmathsetmacro{\bottomLabel}{int(\top - \bottom + 1)}  
            \node[below] at (\i, \minbottom - 2) {\bottomLabel};  
        \fi
    \fi
    

    \draw[thick] (\i, \bottom) -- (\i, \top);
    
    \pgfmathsetmacro{\n}{\top-\bottom} 
    
    \foreach \k in {0,...,\n} {
        \pgfmathsetmacro{\y}{\bottom + \k}
        \filldraw (\i, \y) circle (1.75pt); 
    }

    \ifnum\tail=1
        \draw[thick] (\i, \minbottom) -- (\i, \bottom);
        \pgfmathsetmacro{\m}{\bottom-\minbottom - 1}
        \foreach \k in {0,...,\m} {
            \pgfmathsetmacro{\y}{\minbottom + \k}
            \filldraw (\i, \y) circle (1.75pt); 
        }
        \draw[thick] (\i, \minbottom) -- (\i, \minbottom - 0.5); 
        
        \node at (\i, \minbottom - 0.65) {\scalebox{1}{\vdots}};
    \fi

    \ifnum\tip=1
        \draw[thick, -stealth] (\i, \top) -- ++(0,0.7); 
    \fi

}



\end{tikzpicture}

\caption{Some $w_1$-periodic $h_1$-submodules in the classical Adams $E_\infty$-page.  Numbers along the top are $(\mathrm{coweight}+ 1)$. Numbers along the bottom are $h_1$-exponents.}
\label{fig:w1-periodic}
\end{figure}
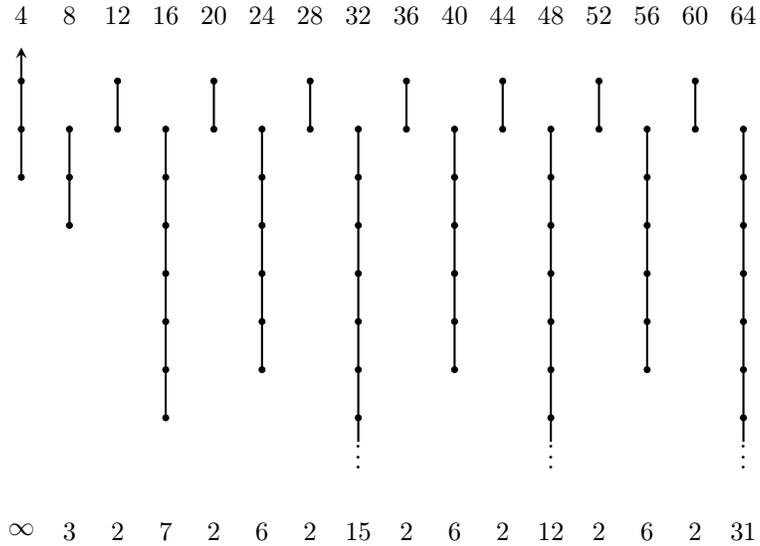

\cref{fig:w1-periodic} displays the contents of \cref{thm:C-homotopy} and \cref{thm:C-homotopy-beta} graphically.  The figure is qualitatively similar to \cref{fig:w1-periodic-algebraic}, but the $h_1$-coexponents of the $h_1$-submodules take slightly different values.  The numbers along the top of the figure indicate the coweight plus one.  The numbers along the bottom of the figure indicate $h_1$-coexponents.  As in \cref{fig:w1-periodic-algebraic},  the elements are arranged such that every horizontal row is a $w_1$-periodic family.  The element $\frac{\overline{w_1^{4k}}}{h_0}$ lies at the top of column $8k$, and the element $\overline{w_1^{4k+2}}$ lies at the top of column $8k+4$.  Beware that the vertical alignment of \cref{fig:w1-periodic} does not reflect the Adams filtration.  On a standard Adams chart, $w_1$-periodic families lie along lines of slope $\frac{1}{5}$, not along horizontal lines.  As in \cref{fig:w1-periodic-algebraic}, the vertical alignments of the columns indicate why the coweights congruent to 7 modulo 8 must be studied somewhat separately from the coweights congruent to 3 modulo 8.

\section{The \texorpdfstring{$\C$}{C}-motivic Burklund--Xu spectral sequence}
\label{sctn:BXss}

Having stated the main theorems of this manuscript, we now start on the technical work required to compute the coexponents that are given in \cref{thm:C-homotopy} and \cref{thm:C-homotopy-beta}.

In this section, we describe the $\C$-motivic version of the Burklund--Xu spectral sequence.  The classical version of this spectral sequence is described in detail in \cite{BX23} and \cite{Benson24}.

The $E_1$-page of the $\C$-motivic Burklund--Xu spectral sequence takes the form
\begin{equation}
\label{eq:BX-E1}
\F_2[\tau][q_0, q_1, \ldots] \otimes_{\F_2} H^{***}\cB.
\end{equation}
Recall from \cref{sctn:background} that the dual $\cB_{**}$ of $\cB$ is the  subalgebra 
\[
\cB_{**} =  \F_2[\xi_1, \xi_2, \ldots ]
\]
of the dual $\C$-motivic Steenrod algebra $\cA_{**}$.  Also, \cref{Ddi} tells us that $H^{***} \cB$ is isomorphic to $H^{**} \cA^{\cl}$, but beware that this isomorphism shifts degrees as in \cref{eq:Ddi}.

\begin{rem}
Note that $\cB$ is similar to the construction $\cP^{\mathrm{mot}}$ of \cite[Construction 2.7]{BX23}.  The difference is that we do not include $\tau$ in $\cB$.
\end{rem}

The $\C$-motivic Burklund--Xu spectral sequence converges to the $E_2$-page of the $\C$-motivic Cartan--Eilenberg spectral sequence, which further converges to $H^{***} \cA$. However, the $\C$-motivic and classical versions of the Burklund--Xu spectral sequence are essentially identical.  After a careful inspection of degrees, one finds that the $\C$-motivic Burklund--Xu differentials involve no powers of $\tau$.  In effect, the $\C$-motivic Burklund--Xu spectral sequence is merely the classical Burklund--Xu spectral sequence tensored over $\F_2$ with $\F_2[\tau]$.

In this manuscript, we are only interested in a small portion of the $\C$-motivic Burklund--Xu spectral sequence.  We will study only the part of the spectral sequence in Chow degrees zero and one.  As observed in \cite[Proposition 10]{Benson24}, the $\C$-motivic Burklund--Xu spectral sequence splits as a direct sum of its subobjects in each fixed Chow degree.  In terms of \cref{eq:BX-E1}, the Chow degree is determined by the facts that each $q_i$ has Chow degree $1$; $\tau$ has Chow degree $2$; and $H^{***}\cB$ is concentrated in Chow degree zero.

\subsection{The \texorpdfstring{$\C$}{C}-motivic Burklund--Xu spectral sequence in Chow degree zero}

First we study the Chow degree zero part of the $\C$-motivic Burklund--Xu spectral sequence. This part of the $E_1$-page is equal to
\[
1^{BX} \cdot H^{***}\cB.
\]
We use the notation $1^{BX} \cdot (-)$ to emphasize that we are considering elements of the Burklund--Xu $E_1$-page.  This helps to disambiguate between multiple interpretations of the same symbol.

\begin{rem}
\label{rem:dr-linear}
Recall from \cite[Proposition 15]{Benson24} that the Burklund--Xu differentials are $H^{***}\cB$-linear, in the sense that
\begin{equation}
\label{eq:dr-linear}
d_r(y \cdot x) = d_r(y) \cdot x
\end{equation}
for all $y$ in the $E_r$-page and all $x$ in $H^{***}\cB$.
\end{rem}

\cref{Ddi} was originally proved in \cite[Theorem 2.19]{Isa19}. We provide an alternative proof using the Burklund--Xu spectral sequence. This new proof is arguably more complicated than the simple observation that $\cA^{\cl}$ and $\cB$ are isomorphic as Hopf algebras.  We offer this proof for its consistency with the techniques used elsewhere in this manuscript.

\begin{proof}[Proof of \cref{Ddi}]
\label{proof:Cd0}
One consequence of \cref{eq:dr-linear} is that elements of the form $1^{BX} \cdot x$ are permanent cycles in the Burklund--Xu spectral sequence \cite[Proposition 15]{Benson24}.  Therefore, the spectral sequence collapses in Chow degree zero.

Moreover, the $\C$-motivic Cartan--Eilenberg spectral sequence in Chow degree zero also collapses.  There are two main points to this observation.  First, $\C$-motivic Cartan--Eilenberg differentials preserve Chow degree, so we need only consider differentials whose sources and targets are both in Chow degree zero.  Second, $\C$-motivic Cartan--Eilenberg differentials take values in multiples of $\tau$ \cite[Theorem 2.8]{BX23}, and there are no multiples of $\tau$ in Chow degree zero.

Hence, the $\C$-motivic Burklund--Xu $E_1$-page in Chow degree zero is isomorphic to the $\C$-motivic Cartan--Eilenberg $E_2$-page in Chow degree zero, which is isomorphic to the part of $H^{***} \cA$ in Chow degree zero.
\end{proof}

\subsection{The \texorpdfstring{$\C$}{C}-motivic Burklund--Xu spectral sequence in Chow degree one}
\label{sctn:BX-Cd1}

We now study the $\C$-motivic Burklund--Xu spectral sequence in Chow degree one in more detail.

\begin{prop}
\label{prop:Cd1}
There is an isomorphism from the $\C$-motivic Cartan--Eilenberg $E_2$-page in Chow degree one to the part of $H^{***} \cA$ in Chow degree one.
\end{prop}

\begin{proof}
As in the proof of \cref{Ddi} just above, $\C$-motivic Cartan--Eilenberg differentials take values in multiples of $\tau$.  Therefore, the values of such differentials have Chow degree of at least $2$.
\end{proof}

After restricting to Chow degree one, the $\C$-motivic Burklund--Xu spectral sequence takes the form
\begin{equation}\label{AC1}
    \bigoplus_{i\geq 0} q_i \cdot H^{***}\cB  \quad\Longrightarrow\quad H^{***}\cA \big|_{s+f-2w=1}.
\end{equation}
The $E_1$-page consists of a sum of copies of $H^{***}\cB$.  Because of \cref{prop:Cd1}, the spectral sequence converges all the way to the part of $H^{***}\cA$ that lies in Chow degree one.

Because of \cref{Ddi}, we find it useful to adopt classical notation for the cohomology of $\cB$.  This can be confusing at first, but we find that it helps us carry out the computations and simplifies formulas in the long run.  Our spectral sequence now takes the form
\begin{equation}\label{AC1.1}
\bigoplus_{i\geq 0}q_i \cdot H^{**}\cA^{\cl} \quad\Longrightarrow\quad H^{***}\cA \big|_{s+f-2w=1}.
\end{equation}

Moreover, we change the gradings on the spectral sequence to be compatible with this use of classical notation.  In the usual motivic grading, our spectral sequence only has elements in degrees of the form $\left( s, f, \frac{s+f-1}{2} \right)$ because it is concentrated in Chow degree one.  Such elements are now placed in degree $\left( \frac{s - f + 1}{2}, f \right)$ in our modified grading.

This change of grading has the following practical effects.  First, the elements $q_i$ have degree $(2^{i+1} - 2, 1, 2^i - 1)$ in the traditional motivic grading.  In our new grading, these elements now have degree $(2^i - 1, 1)$.  Second, elements of $H^{**}\cA^{\cl}$ now have their familiar degrees.

\begin{ex}
The element $q_2 \cdot h_0^2$ has degree $(3,3) = (3,1) + (0,2)$.  This degree converts back to $(8,3,5)$ in the traditional motivic grading, and in fact the element $q_2 \cdot h_0^2$ represents $c_0$.
\end{ex}

In our regrading, all differentials have degree $(-1,1)$, since we are not explicitly including the Burklund--Xu filtration in the grading.

Recall also from \cite[Proposition 14]{Benson24} \cite[Lemma 4.10]{BX23} that there are $d_1$ differentials
\begin{equation}
\label{eq:BX-d1}
\begin{split}
& d_1(q_i)=q_{i-1} \cdot h_{i-1} \quad \mathrm{ if } \quad i > 0. \\
& d_1(q_0) = 0.
\end{split}
\end{equation}
These formulas, together with $H^{**}\cA^{\cl}$-linearity described in \cref{rem:dr-linear}, completely determine the Burklund--Xu $d_1$ differential.

\begin{rem}
\label{rem:higher-diff}
There is also a formula involving threefold Massey products for the $d_2$ differential \cite[Proposition 14]{Benson24} \cite[Lemma 4.10]{BX23}.  For $i \geq 2$ we have
\[
d_2(q_i \cdot x) = q_{i-2} \cdot \langle h_{i-2}, h_{i-1}, x \rangle,
\]
where the Massey product occurs in $H^{**}\cA^{\cl}$.  The well-definedness of the bracket is assured by the assumption that $d_1(q_i \cdot x) = 0$, and the indeterminacy is zero in the $E_2$-page because it lies in the image of the $d_1$ differential.

Similarly, there are formulas involving higher order Massey products for the higher differentials.  Specifically, for $i \geq r$ we have
\[
d_r(q_i \cdot x) = q_{i-r} \cdot \langle h_{i-r}, \ldots, h_{i-2}, h_{i-1}, x \rangle.
\]
See \cite[Lemma 4.10]{BX23} for an explicit statement of the formula when $r = 3$.  The earlier differentials take care of the well-definedness and indeterminacy of the bracket.  These formulas are algebraic analogues of the homotopical Toda bracket formulas of \cite{CF17}.  We do not discuss the proofs of these higher differential formulas because we do not need them.
\end{rem}

\section{Analysis of the Burklund--Xu spectral sequence in a range}
\label{sctn:BXss-analysis}

The goal of this section is to study the Burklund--Xu spectral sequence explicitly in a range of degrees.  Our analysis will lead to a family of non-zero products in the cohomology of the $\C$-motivic Steenrod algebra in Chow degree one.

Recall from \cref{sctn:background} that $s-2f$ is the $v_1$-intercept of an element in stem $s$ and Adams filtration $f$.

\begin{ex}
\label{ex:vn-intercept}
The element $q_n$ of the Burklund--Xu spectral sequence has degree $(2^n-1,1)$, so its $v_1$-intercept is $2^n-3$.
\end{ex}

\begin{rem}
\label{rem:intercept-multiplicative}
The $v_1$-intercept is multiplicative in the sense that the $v_1$-intercept of $xy$ is equal to the sum of the $v_1$-intercepts of $x$ and $y$.
\end{rem}

\subsection{The cohomology of the classical Steenrod algebra near the vanishing line}

We restate a classical result of Adams about the vanishing of the Adams $E_2$-page above a line of slope $\frac{1}{2}$.

\begin{prop}
\label{prop:Adams-vanishing}
\cite{Ada66}
The cohomology of the classical Steenrod algebra is concentrated in degrees with $v_1$-intercept greater than or equal to $-3$, with the exception of the elements $h_0^k$ in degree $(0,k)$ for $k \geq 2$.
\end{prop}

On a standard Adams chart, \cref{prop:Adams-vanishing} says that the chart is concentrated below a line of slope $\frac{1}{2}$ and $x$-intercept $-3$, with the exception of the elements $h_0^k$ for $k \geq 2$.

\begin{lem}
\label{prop:low-v1-intercept}
\cref{table:low-v1-intercept} lists all elements in the cohomology of the classical Steenrod algebra whose $v_1$-intercepts equal $-2$, $0$, or $1$.
\end{lem}

\begin{proof}
The degrees under consideration lie within the region in which the Adams periodicity operator $P$ of degree $(8,4)$ is an isomorphism \cite{Ada66} \cite{Li20}.  The result follows from explicit low-dimensional computations in degrees up to $(12,6)$.
\end{proof}

\begin{longtable}{rll}
\caption{Elements in the cohomology of the classical Steenrod algebra with $v_1$-intercepts equal to $-2$, $0$, and $1$
\label{table:low-v1-intercept}
} \\
\toprule
$v_1$-intercept & $(s,f)$ & element \\
\midrule \endfirsthead
\caption[]{Elements in the cohomology of the classical Steenrod algebra with $v_1$-intercepts equal to $-2$, $0$, and $1$}\\
\toprule
$v_1$-intercept & $(s,f)$ & element \\
\midrule \endhead
\bottomrule \endfoot
$-2$ & $(0,1)$ & $h_0$ \\
$-2$ & $(8k + 2, 4k + 2)$ & $P^k h_1^2$ \\
$0$ & $(0,0)$ & $1$ \\
$1$ & $(8k+3, 4k+1)$ & $P^k h_2$ \\
$1$ & $(8k+7, 4k+3)$ & $P^k h_0^2 h_3$ \\
$1$ & $(8k+9, 4k+4)$ & $P^k h_1 c_0$
\end{longtable}

\subsection{The Burklund--Xu \texorpdfstring{$E_1$}{E1}-page in certain degrees}

We need some technical results that classify the elements in the Burklund--Xu $E_1$-page in various degrees.

\begin{longtable}{lll}
\caption{Elements of the Burklund--Xu $E_1$-page in Chow degree one with $v_1$-intercept equal to $-1$
\label{table:BX-E1}
} \\
\toprule
$(s,f)$ & element & condition\\
\midrule \endfirsthead
\caption[]{Some elements of the Burklund--Xu $E_1$-page in Chow degree one with $v_1$-intercept equal to $-1$}\\
\toprule
$(s,f)$ & element & condition\\
\midrule \endhead
\bottomrule \endfoot
$(8k+3, 4k+2)$ & $q_0 \cdot P^k h_2$ & $k \geq 0$ \\
$(8k+7, 4k+4)$ & $q_0 \cdot P^k h_0^2 h_3$ & $k \geq 0$ \\
$(8k+9, 4k+5)$ & $q_0 \cdot P^k h_1 c_0$ & $k \geq 0$ \\
$(1,1)$ & $q_1 \cdot 1$ \\
$(8k + 5, 4k + 3)$ & $q_2 \cdot P^k h_1^2$ & $k \geq 0$ \\
$(2^n-1, 2^{n-1})$ & $q_n \cdot h_0^{2^{n-1} - 1}$ & $n \geq 2$ \\
\end{longtable}

\begin{lem}
\label{prop:v1-intercept--1}
\cref{table:BX-E1} lists all elements of the Burklund--Xu $E_1$-page in Chow degree one with $v_1$-intercept equal to $-1$.
\end{lem}

\begin{proof}
We are searching for expressions of the form $q_n \cdot x$ such that $x$ belongs to $H^{**} \cA^{\cl}$, and $q_n \cdot x$ has $v_1$-intercept $-1$.  Since $q_n$ has $v_1$-intercept $2^n-3$, our desired $x$ must have $v_1$-intercept $2 - 2^n$.

If $n \geq 3$, then $ 2 - 2^n \leq -6$.  By \cref{prop:Adams-vanishing}, $x$ must be zero or of the form $h_0^j$.  More specifically, the elements $q_n \cdot h_0^{2^{n-1} - 1}$ are the only elements of the form $q_n \cdot x$ such that $n \geq 3$, the $v_1$-intercept is $-1$, and the Chow degree is $1$.

Next, consider elements of the form $q_2 \cdot x$.  Then $x$ has $v_1$-intercept $-2$ since $q_2$ has $v_1$-intercept $1$.  By \cref{prop:low-v1-intercept}, $x$ must be $h_0$ or of the form $P^k h_1^2$ for some $k \geq 0$.

Now consider elements of the form $q_1 \cdot x$.  Then $x$ has $v_1$-intercept $0$ since $q_1$ has $v_1$-intercept $-1$.  By \cref{prop:low-v1-intercept}, $x$ must be $1$.

Finally, consider elements of the form $q_0 \cdot x$.  Then $x$ has $v_1$-intercept $1$ since $q_0$ has $v_1$-intercept $-2$.  By \cref{prop:low-v1-intercept}, $x$ must be of the form $P^k h_2$, $P^k h_0^2 h_3$, or $P^k h_1 c_0$ for some $k \geq 0$.
\end{proof}

\begin{lem}
\label{lem:v1-intercept--7}
Consider the Burklund--Xu $E_1$-page in Chow degree one and $v_1$-intercept $-7$.  Under these restrictions, the only elements are $q_n \cdot h_0^{2^{n-1} + 2}$ in degree $(2^n - 1, 2^{n-1} + 3)$ for $n \geq 1$.
\end{lem}

\begin{proof}
We are searching for expressions of the form $q_n \cdot x$ with $x$ in $H^{**}\cA^{\cl}$, such that $q_n \cdot x$ has $v_1$-intercept $-7$.  Since $q_n$ has $v_1$-intercept $2^n-3$, our desired $x$ must have $v_1$-intercept $-4 - 2^n$.

If $n \geq 1$, then $-4-2^n \leq -6$.  By \cref{prop:Adams-vanishing}, $x$ can only be a power of $h_0$.
\end{proof}

\subsection{The elements \texorpdfstring{$h_0 h_2 \cdot h_2 g^k$}{h0h2h2gk}}

We now study a $w_1$-periodic family of products of the form $h_0 h_2 \cdot h_2 g^k$ for $k \geq 1$ in the cohomology of the $\C$-motivic Steenrod algebra.  We will show that all of these products are non-zero.  Note that the Chow degree of $h_0h_2\cdot h_2g^k$ is $1$.  We will identify their representatives in the Burklund--Xu spectral sequence, and we will show that these representatives are non-zero permanent cycles.

\begin{lem}
\label{prop:P^kh1^2.v0}
For $k \geq 1$, the elements $q_0 \cdot P^k h_1^2$ in degree $(8k+2, 4k+3)$ are non-zero permanent cycles in the Burklund--Xu spectral sequence.
\end{lem}

\begin{proof}
The elements $q_0 \cdot P^k h_1^2$ are in the lowest possible Burklund--Xu filtration, and the Burklund--Xu $d_r$ differential decreases the filtration by $r$.  Therefore, these elements are permanent cycles.

To prove that $q_0 \cdot P^k h_1^2$ cannot be the target of a differential, we need to rule out all possible differentials that might hit those elements.  The degree of $q_0 \cdot P^k h_1^2$ is $(8k+2, 4k+3)$, and the sources of possible differentials have degree $(8k+3, 4k+2)$.  The $v_1$-intercept of such elements is $-1$.  \cref{prop:v1-intercept--1} classifies elements with $v_1$-intercept $-1$.  Of those elements, only the elements $q_0 \cdot P^k h_2$ lie in a degree of the form $(8k+3,4k+2)$ with $k \geq 1$.  The Burklund--Xu filtration of $q_0 \cdot P^k h_2$ is $0$, so it cannot support a differential.
\end{proof}

\begin{rem}
\label{rem:P^kh1^2.v0}
The condition $k \geq 1$ in \cref{prop:P^kh1^2.v0} is necessary.  If $k = 0$, then there is an element $q_2 \cdot h_0$ in degree $(3,2) = (8k+3, 4k+2)$.  In fact, there is a differential
\[
d_2(q_2 \cdot h_0) = q_0 \cdot h_1^2.
\]
As in \cref{rem:higher-diff}, this differential is equivalent to the Massey product
\[
h_1^2 = \langle h_0, h_1, h_0 \rangle.
\]
Therefore, $q_0 \cdot h_1^2$ is zero in the $E_\infty$-page.
\end{rem}

\begin{prop}
\label{prop:h0h2.h2g^k}
For all $k \geq 1$, the product $h_0 h_2 \cdot h_2 g^k$ in degree $(20k + 6, 4k + 3)$ is non-zero in $H^{***} \cA$.
\end{prop}

\begin{proof}
The proof consists in converting the non-zero permanent cycles $q_0 \cdot P^k h_1^2$ of \cref{prop:P^kh1^2.v0} from classical notation into motivic notation.  The factor $P^k h_1^2$ corresponds to $h_2 \cdot h_2 g^k$, while the factor $q_0$ corresponds to ($\C$-motivic) $h_0$.
\end{proof}

\begin{rem}
\label{rem:h0h2.h2g^k}
The condition $k \geq 1$ in \cref{prop:h0h2.h2g^k} is necessary.  When $k = 0$, the product $h_0 h_2 \cdot h_2$ is in fact zero.  This corresponds to the differential discussed in \cref{rem:P^kh1^2.v0}.
\end{rem}

\subsection{Further structure of the Burklund--Xu spectral sequence in
Chow degree one}
\label{subsctn:BX-further}

For the interested reader, we provide some details about the full Burklund--Xu spectral sequence for Chow degree one in the range through coweight $30$.  The main results of this manuscript in no way depend on the material in this subsection, and we do not offer complete proofs for all claims.

\begin{figure}[H]
\begin{center}
\includegraphics[trim={0cm, 0cm, 0cm, 0cm},clip,page=1,scale=0.44]{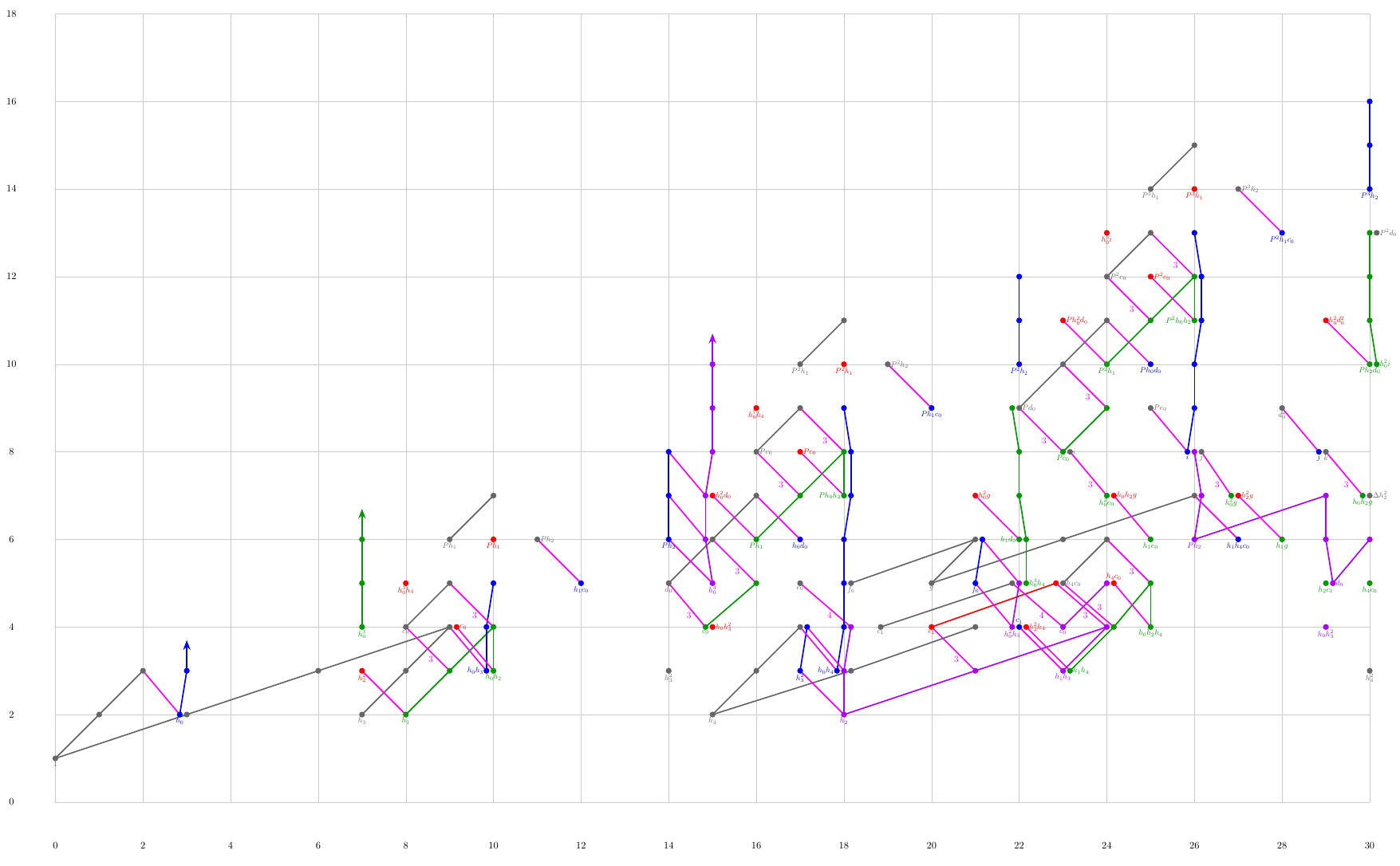}
\caption{The $E_2$-page of the Burklund--Xu spectral sequence for the cohomology of the $\C$-motivic Steenrod algebra in Chow degree one.
\label{fig:A-E2}
}
\end{center}
\end{figure}

\cref{fig:A-E2} displays the Burklund--Xu $E_2$-page in Chow degree one up to coweight $30$.  Dots represent copies of $\mathbb{F}_2$.  The labeled elements, together with \cref{tab:colors}, indicate the names of elements in the form $q_n \cdot x$.  The value of $n$ is given in \cref{tab:colors}, and the value of $x$ is displayed in the label.

\begin{longtable}{ll}
\caption{Color coding for \cref{fig:A-E2}
\label{tab:colors}
} \\
\toprule
color & generator \\
\midrule \endfirsthead
\caption[]{Color coding for \cref{fig:A-E2}}\\
\toprule
color & generator \\
\midrule \endhead
\bottomrule \endfoot
{\color{gray} gray} & $q_0$ \\
{\color{red} red} & $q_1$ \\
{\color{blue} blue} & $q_2$ \\
{\color{green} green} & $q_3$ \\
{\color{purple} purple} & $q_4$ 
\end{longtable}

As usual, vertical lines indicate multiplications by $h_0$, lines of slope $1$ indicate multiplications by $h_1$, and lines of slope $\frac{1}{3}$ indicate multiplications by $h_2$.  Vertical arrows indicate infinitely many $h_0$ multiplications.  Because of our notational conventions, these lines represent multiplications by $h_1$, $h_2$, and $h_3$ respectively in the cohomology of the $\C$-motivic Steenrod algebra.

The higher Burklund--Xu $d_r$ differentials for $r \geq 2$ are displayed in \cref{fig:A-E2} as magenta lines of slope $-1$.  All differentials have the same slope because the Burklund--Xu filtration is suppressed in our charts.  An unlabeled magenta line indicates a $d_2$ differential, while a magenta line labelled $r$ indicates a $d_r$ differential for $r \geq 3$.

In terms of colors, note that the $d_r$ differential moves $r$ stages up in \cref{tab:colors}.  In particular, in this range the spectral sequence converges at the $E_5$-page.

\begin{figure}[H]
\begin{center}
\includegraphics[trim={0cm, 0cm, 0cm, 0cm},clip,page=1,scale=0.44]{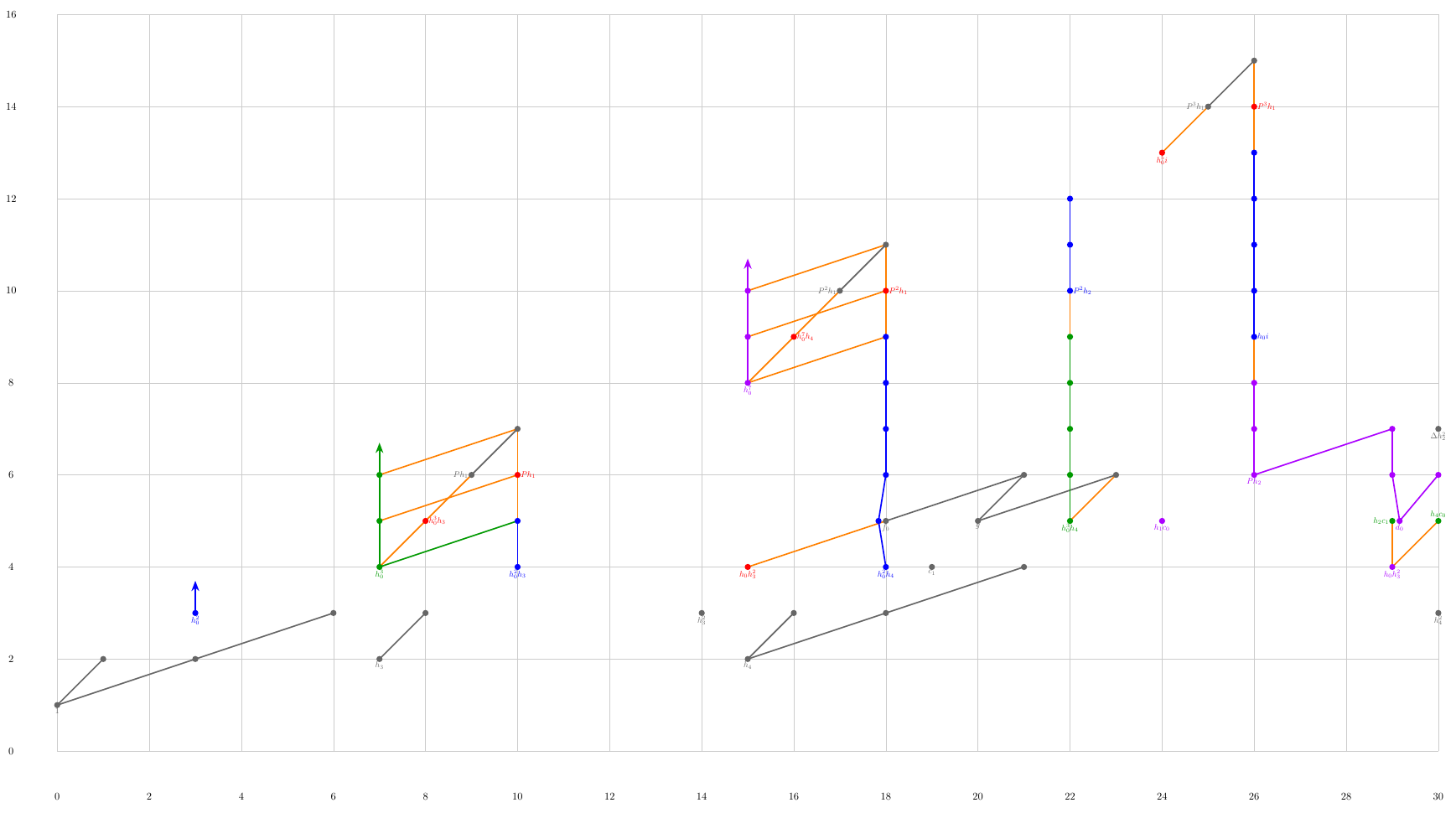}
\caption{The $E_\infty$-page of the Burklund--Xu spectral sequence for the cohomology of the $\C$-motivic Steenrod algebra in Chow degree one
\label{fig:A-Einfty}
}
\end{center}
\end{figure}

\cref{fig:A-Einfty} displays the Burklund--Xu $E_\infty$-page in Chow degree one up to coweight 30.  The interpretation of dots, colors, labels, and lines is the same as in \cref{fig:A-E2}.  Orange lines indicate extensions in the cohomology of the $\C$-motivic Steenrod algebra that are hidden in the Burklund--Xu $E_\infty$-page.

\begin{landscape}

\begin{figure}[H]
\begin{center}
\includegraphics[trim={0cm, 0cm, 0cm, 0cm},clip,page=1,scale=0.28]{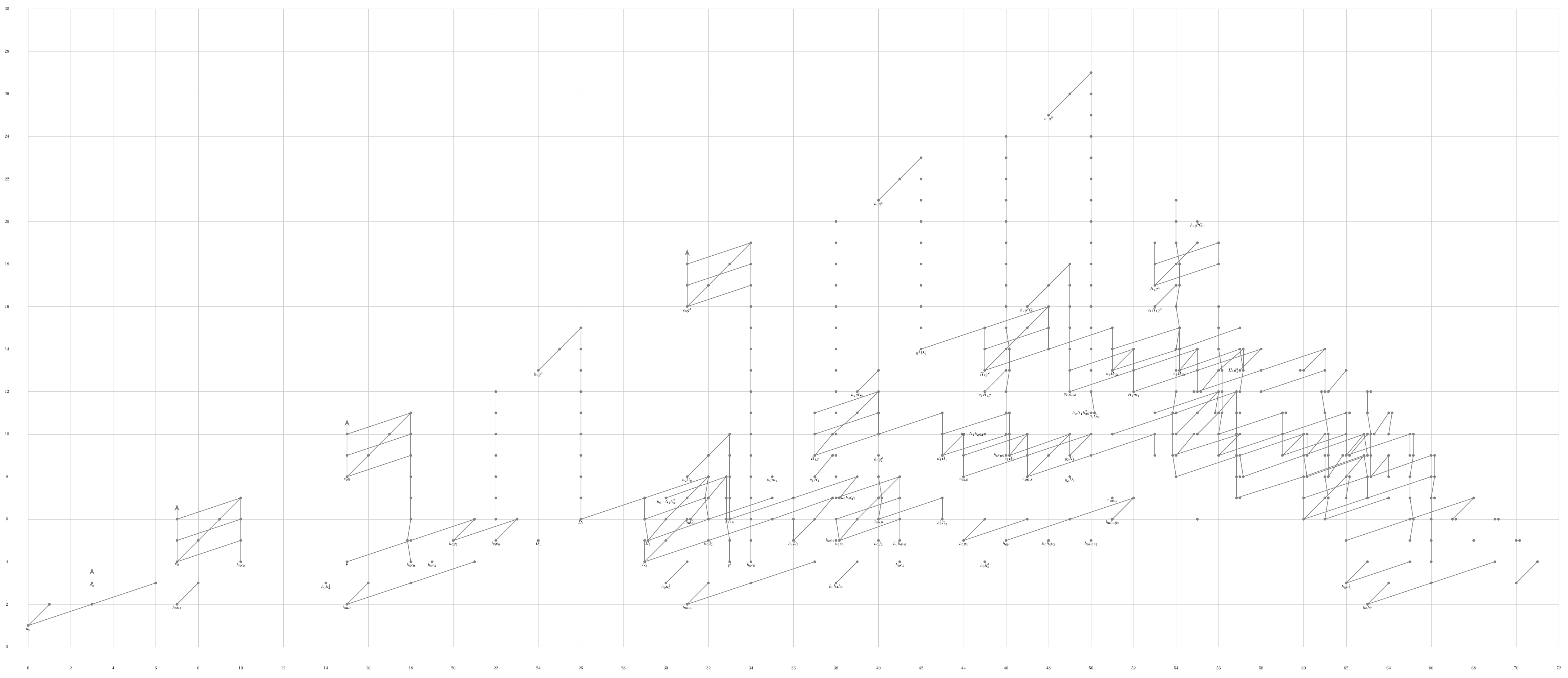}
\caption{The cohomology of the $\C$-motivic Steenrod algebra in Chow degree one
through the 140-stem
\label{fig:Guozhen}
}
\end{center}
\end{figure}

\end{landscape}

Beyond coweight 30, we have not studied the structure of the spectral sequence closely.  However, using raw data through the 140-stem \cite{IWX22b}, we extracted a complete description of the part of the cohomology of the $\C$-motivic Steenrod algebra in Chow degree one in that same range.  The results are displayed in \cref{fig:Guozhen}.  The periodic structure along the top of the chart along a line of slope $\frac{1}{2}$ is the primary focus of this manuscript.

\section{Adams differentials}
\label{sctn:Adams-differentials}

The goal of this section is to establish an infinite family of non-zero $\C$-motivic Adams differentials, as well as to establish several infinite families of permanent cycles in the $\C$-motivic Adams spectral sequence.  \cref{tab:elements} lists some of the specific elements that we will study.  All of these elements have Chow degree zero.  The third column of the table describes the elements in terms of the $\C$-motivic Burklund--Xu spectral sequence.

\begin{longtable}{lll}
\caption{Some $w_1$-periodic elements in the cohomology of the $\C$-motivic Steenrod algebra
\label{tab:elements}
} \\
\toprule
degree & element & definition \\
\midrule \endfirsthead
\caption[]{Some $w_1$-periodic elements in the cohomology of the $\C$-motivic Steenrod algebra}\\
\toprule
degree & element & definition \\
\midrule \endhead
\bottomrule \endfoot
$(20k+7, 4k+1, 12k+4)$ & $h_3 g^k$ & $1^{BX} \cdot P^k h_2$  \\
$(20k+3, 4k+1, 12k+2)$ & $h_2 g^k$ & $1^{BX} \cdot P^k h_1$  \\
$(40k+53, 8k+7, 24k+30)$ & $i_1 g^{2k}$ & $1^{BX} \cdot P^{2k} i$  \\
$(80k+107, 16k+13, 48k+60)$ & $\Delta_1 h_1^2 i_1 g^{4k}$ & $1^{BX} \cdot P^{4k} \Delta h_0^2 i$ \\
\end{longtable}

\subsection{Some Massey products in the cohomology of the \texorpdfstring{$\C$}{C}-motivic Steenrod algebra}

In order to analyze Adams differentials on specific elements, we will use some Massey product decompositions in the cohomology of the $\C$-motivic Steenrod algebra.

\begin{prop}
\label{prop:Massey}
\cref{tab:Massey} shows some infinite families of non-zero Massey products in the cohomology of the $\C$-motivic Steenrod algebra.  All of the indeterminacies are zero.
\end{prop}

\begin{longtable}{llll}
\caption{Some Massey products in the cohomology of the $\C$-motivic Steenrod algebra
\label{tab:Massey}
} \\
\toprule
degree & element & Massey product & condition \\
\midrule \endfirsthead
\caption[]{Massey products in the cohomology of the $\C$-motivic Steenrod algebra}\\
\toprule
degree & element & Massey product & condition \\
\midrule \endhead
\bottomrule \endfoot
$(20k+67, 4k+13, 12k+40)$ & $h_3 g^{k+3}$ & $\langle h_3 g^k, h_1^3, h_1^3 i_1 \rangle$ & $k \geq 0$ \\
$(20k+63, 4k+13, 12k+38)$ & $h_2 g^{k+3}$ & $\langle h_2 g^k, h_1^3, h_1^3 i_1 \rangle$ & $k \geq 0$ \\
$(20k+28, 4k+6, 12k+17)$ & $h_1 h_3 g^{k+1}$ & $\langle h_1 h_3 g^k, h_1^3, h_1 h_4 \rangle$ & $k \geq 0$ \\
$(40k+93, 8k+15, 24k+54)$ & $i_1 g^{2k+2}$ & $\langle i_1 g^{2k}, h_1^7, h_1 h_5 \rangle$ & $k \geq 0$ \\
$(80k \mathord{+} 187, 16k \mathord{+} 29, 48k \mathord{+} 108)$ & $\Delta_1 h_1^2 i_1 g^{4k+4}$ & $\langle \Delta_1 h_1^2 i_1 g^{4k}, h_1^{15}, h_1 h_6 \rangle$ & $k \geq 0$ \\
$(20k+66, 4k+15, 12k+40)$ & $h_0 h_2 \cdot h_2 g^{k+3}$ & $\langle h_0 h_2 \cdot h_2 g^k, h_1^3, h_1^3 i_1 \rangle$ & $k \geq 1$ \\
\end{longtable}

The Massey products in \cref{tab:Massey} are recursive, in the sense that each member of each infinite family equals a Massey product that involves a previous member of that same family.

\begin{proof}
For the first five families, all of the relevant elements have Chow degree zero.  The Massey products are analogous to the classical non-zero Massey products
\[
\begin{split}
P^{k+3} h_2 = \langle P^k h_2, h_0^3, h_0^3 i \rangle \\
P^{k+3} h_1 = \langle P^k h_1, h_0^3, h_0^3 i \rangle \\
P^{k+1} h_0 h_2 = \langle P^{k} h_0h_2 ,h_0^3,h_0h_3\rangle \\
P^{2k+2} i = \langle P^{2k} i, h_0^7, h_0 h_4\rangle \\
P^{4k+4} \Delta h_0^2 i = \langle P^{4k} \Delta h_0^2 i, h_0^{15}, h_0 h_5 \rangle \\
\end{split}
\]
under the isomorphism of \cref{Ddi} from  $H^{**} \cA^{\cl}$ to the Chow degree zero part of $H^{***} \cA$.

For the last family, we analyze the indeterminacy first.  The indeterminacy of the bracket consists of linear combinations of expressions of the form $h_0 h_2 \cdot h_2 g^k \cdot x$, where $x$ has degree $(60,12,36)$; and expressions of the form $y \cdot h_1^3 i_1$, where $y$ has degree $(20k+10, 4k+5, 12k+7)$.  There are no non-zero elements in degree $(60,12,36)$.

We now consider elements of degree $(20k+10, 4k+5, 12k+7)$.  Note that the Chow degree of such elements is $1$, so our earlier analysis of the Burklund--Xu spectral sequence in \cref{sctn:BXss-analysis} is relevant.  After converting to the degrees that we used in that spectral sequence, we are looking for elements of degree $(8k+3, 4k+5)$.  The $v_1$-intercept of such elements is $-7$.  \cref{lem:v1-intercept--7} classifies elements of the Burklund--Xu $E_1$-page with $v_1$-intercept $-7$, but none of those elements has a degree of the form $(8k+3, 4k+5)$ with $k \geq 1$.  Therefore, there are no non-zero elements of degree $(20k+10, 4k+5, 12k+7)$, and the Massey product has no indeterminacy.

Finally, we must determine the unique element that is contained in the Massey product. We have
\[
\langle h_0 h_2 \cdot h_2 g^k, h_1^3, h_1^3 i_1 \rangle 
= h_0 h_2 \langle h_2 g^k, h_1^3, h_1^3 i_1 \rangle
= h_0 h_2 \cdot h_2 g^{k+3},
\]
where the first equality holds because there is no indeterminacy, and the last equality is the second family of Massey products in \cref{tab:Massey}.  \cref{prop:h0h2.h2g^k} shows that the Massey product is non-zero.
\end{proof}

\begin{rem}
The $k=0$ case of the last family in \cref{tab:Massey} behaves differently.  In that case, we have the Massey product $\langle 0, h_1^3, h_1^3 i_1 \rangle$ in degree $(66,15,40)$ because $h_0 h_2^2 = 0$.  In that degree, there is a relation 
\[
h_0 h_2 \cdot h_2 g^3 = h_1^3 i_1 \cdot h_1^2 c_0.
\]
So the non-zero element $h_0 h_2 \cdot h_2 g^3$ does belong to the bracket, but the indeterminacy makes this fact not useful.
\end{rem}

\subsection{Adams differentials \texorpdfstring{$d_2(h_3 g^k)$}{d2(h3gk)}}

The goal of this section is to prove that the classes $h_3g^k$, for $k \geq 1$, in the $\C$-motivic Adams $E_2$-page support Adams differentials.  Explicit computation \cite{IWX20} shows that there are Adams differentials $d_2(h_3 g^k) = h_0 h_2 \cdot h_2 g^k$ for $k \leq 5$.

\begin{prop}\label{Cd2}
For all $k \geq 1$, there is a non-zero $\C$-motivic Adams $d_2$ differential $$d_2(h_3g^k)=h_0h_2\cdot h_2g^k.$$ 
\end{prop} 

\begin{proof}
The proof is by induction.  The base cases are $1 \leq k \leq 3$.  These cases follow from explicit low-dimensional computation \cite{IWX20}.

There is a ``higher Leibniz rule'' \cite[Theorem 1.1]{Mos70} that describes how the Adams $d_2$ differential interacts with Massey products in the Adams $E_2$-page.  Applied to the first family of Massey products in \cref{tab:Massey}, the rule says that $d_2(h_3 g^{k+3})$ belongs to
\[
\langle d_2(h_3 g^k), h_1^3, h_1^3 i_1 \rangle +
\langle h_3 g^k, 0, h_1^3 i_1 \rangle +
\langle h_3 g^k, h_1^3, 0 \rangle.
\]
The latter two terms consist of the multiples of $h_3 g^k$ and of $h_1^3 i_1$ in the appropriate degree, and we must determine these multiples.

First, we look for expressions of the form $h_3 g^k \cdot x$, where $x$ has degree $(59,14,36)$.  By inspection, $x$ must be zero.

Second, we look for expressions of the form $y \cdot h_1^3 i_1$, where $y$ has degree $(20k+10, 4k+5, 12k+7)$.  We already considered this situation in the proof of \cref{prop:Massey}, and we found that $y$ must be zero.

We have now determined that $d_2(h_3 g^{k+3})$ belongs to $\langle d_2(h_3 g^k), h_1^3, h_1^3 i_1 \rangle$.  By induction, this bracket equals $\langle h_0 h_2 \cdot h_2 g^k, h_1^3, h_1^3 i_1 \rangle$, which equals $h_0 h_2 \cdot h_2 g^{k+3}$ as shown in \cref{tab:Massey}.
\end{proof}

\begin{rem}
The use of $g^3$-periodicity, rather than $g$-periodicity, is a curious feature of the proof of \cref{Cd2}.  The key point is that there are no elements in $H^{***}\cA$ in degree $(59,14,36)$.  Note that the expression $h_1^2 e_0 g^2$ would have degree $(59,14,36)$ if it existed, but it does not exist in $H^{***}\cA$.  See \cite{KILRZ-e0gk} \cite{Tha21} for information about the non-existence of this element.

If we had used $g$-periodicity, then the indeterminacy in our computations would have involved elements in degree $(19,6,12)$.  In this case, the element $h_1^2 e_0$ is present.  That would result in indeterminacy in the brackets, and the argument would fail.  In fact, the same problem occurs for $g^{2^k}$-periodicity, in which the presence of $h_1^2 e_0 g^{2^k-1}$ creates indeterminacy problems.  We chose to use $g^3$-periodicity because $3$ is the smallest number that is not a power of $2$.
\end{rem}

\subsection{The permanent cycles \texorpdfstring{$h_1h_3g^k$}{h1h3gk}}

Next, we study the elements $h_1 h_3 g^k$, which have Chow degree $0$.  Our goal is to show that they are non-zero permanent cycles in the $\C$-motivic Adams spectral sequence for all $k \geq 0$.

\begin{prop}
\label{h1h3g^k}
The element $h_1 h_3 g^k$ is a non-zero permanent cycle in the $\C$-motivic Adams spectral sequence for all $k \geq 0$.
\end{prop}

\begin{proof}
\cref{prop:no-differentials} shows that $h_1 h_3 g^k$ is non-zero if it does not support any differentials.

The proof is by induction.  The case $k= 0$ is exceptional, but we know from explicit computation that $h_1 h_3$ is a permanent cycle.  The base case $k=1$ says that $h_1 h_3 g$ is a permanent cycle, which we also know from explicit computation \cite{IWX20}.  Let $\alpha_1$ be an element of $\pi_{28,17}$ that is detected by $h_1 h_3 g$.  Note that $\eta^3 \alpha_1$ is zero for degree reasons \cite{IWX22b}.

We inductively define $\alpha_k$ to be an element of the Toda bracket $\langle \alpha_{k-1}, \eta^3, \eta_4 \rangle$.  Part (2) of \cref{thm:C-homotopy-beta} ensures that $\alpha_{k-1} \cdot \eta^2$ is zero, so the bracket is well-defined.

The induction hypothesis is that $\alpha_{k-1}$ is detected by $h_1 h_3 g^{k-1}$.  The Moss convergence theorem \cite[Theorem 1.1]{BK21} implies that the Toda bracket $\langle \alpha_{k-1}, \eta^3, \eta_4 \rangle$ is detected by the Massey product $\langle h_1 h_3 g^{k-1}, h_1^3, h_1 h_4 \rangle$, which equals $h_1 h_3 g^k$ as shown in \cref{tab:Massey}.  The Moss convergence theorem requires that there be no ``crossing differentials''.  In our case, such crossing differentials would have sources with negative Chow degree, but there are no such elements by \cref{prop:Chow-degree-nonneg}.  Therefore, $\alpha_k$ is detected by $h_1 h_3 g^k$.
\end{proof}

\begin{rem}
Extending the notation of the proof, we may define $\alpha_0$ to be an element of $\pi_{8,5}$ that is detected by $h_1 h_3$.  However, $\eta^3 \alpha_0$ is non-zero, and it is detected by $h_1^3 c_0$.  Consequently, the base case for the induction is $\alpha_1$, which is annihilated by $\eta^3$.
\end{rem}

\begin{rem}
The elements $h_1^2 h_3 g^k$ are also non-zero permanent cycles in the $\C$-motivic Adams spectral sequence.  This follows immediately from \cref{h1h3g^k} by multiplication with $h_1$.  This family of elements was previously studied by Andrews \cite{And}.
\end{rem}

\subsection{The permanent cycles \texorpdfstring{$i_1g^{2k}$}{i1g2k}}

\begin{prop}
\label{i1g^2k}
The element $i_1 g^{2k}$ is a non-zero permanent cycle in the $\C$-motivic Adams spectral sequence for all $k \geq 0$.
\end{prop}

\begin{proof}
\cref{prop:no-differentials} shows that $i_1 g^{2k}$ is non-zero if it does not support any differentials.

The proof proceeds by induction.  The base case $k=0$ says that $i_1$ is a permanent cycle, which we know from explicit computation \cite{IWX20}.  Let $\alpha_0$ be an element of $\pi_{53,30}$ that is detected by $i_1$.  Note that $\eta^7 \alpha_0$ is zero for degree reasons \cite{IWX22b}.

We inductively define $\alpha_k$ to be an element of the Toda bracket $\langle \alpha_{k-1}, \eta^7, \eta_5 \rangle$.  Part (2) of \cref{thm:C-homotopy} ensures that $\alpha_{k-1} \cdot \eta^6$ is zero, so the bracket is well-defined.

The induction hypothesis is that $\alpha_{k-1}$ is detected by $i_1 g^{2k-2}$.  The Moss convergence theorem \cite[Theorem 1.1]{BK21} implies that the Toda bracket $\langle \alpha_{k-1}, \eta^7, \eta_5 \rangle$ is detected by the Massey product $\langle i_1 g^{2k-2}, h_1^7, h_1 h_5 \rangle$, which equals $i_1 g^{2k}$ as observed in \cref{tab:Massey}.  The Moss convergence theorem requires that there be no ``crossing differentials''.  In our case, such crossing differentials would have sources with negative Chow degree, but there are no such elements by \cref{prop:Chow-degree-nonneg}.  Therefore, $\alpha_k$ is detected by $i_1 g^{2k}$.
\end{proof}

\begin{rem}
The elements $h_1^m i_1 g^{2k}$ are also non-zero permanent cycles for $0 \leq m \leq 5$.  This follows from \cref{i1g^2k} by multiplication with powers of $h_1$.
\end{rem}

\subsection{The permanent cycle \texorpdfstring{$\Delta_1h_1^2i_1$}{D1h12i1}}
\label{subsctn:D1h12i1}

In this section, we study a particular element $\Delta_1 h_1^2 i_1$ in degree $(107,13,60)$.  Note that this element has Chow degree $0$.  Our goal is to show that it is a permanent cycle.  We already know from \cref{prop:no-differentials} that it cannot be hit by an Adams differential, so our result will show that $\Delta_1 h_1^2 i_1$ detects a non-zero element in $\pi_{107,60}$.

We rely heavily on machine-generated data for the $\C$-motivic Adams $E_2$-page, given in the file \texttt{Adams-motivic-E2-machine.csv} available at \cite{IWX22b}.  We adopt the same notation used in that file, in which $\mathbb{F}_2[\tau]$-module generators of the $\C$-motivic $E_2$-page are named in the form $\texttt{\{a-b\}}$, where $\texttt{a}$ and $\texttt{b}$ are natural numbers.  The value of $\texttt{a}$ indicates the Adams filtration of the element, while the value of $\texttt{b}$ is an arbitrary index whose only purpose is to distinguish different elements.  For example, $\texttt{\{13-300\}}$ refers to the element $\Delta_1 h_1^2 i_1$ under consideration.  In the following lemmas, we will show that $\texttt{\{13-300\}}$ cannot support an Adams differential.

For the convenience of the careful reader, each of the lemmas below is labeled with the specific degree of the element that is the subject of the lemma.

\begin{lem}\label{D1h1^2i1-d2}
$(107, 13, 60)$ 
$d_2 \left( \emph{\texttt{\{13-300\}}} \right) = 0$.
\end{lem}

\begin{proof}
By inspection, the only possible non-zero value for $d_2 \left( \texttt{\{13-300\}} \right)$ is $\texttt{\{15-313\}}$, also known as $h_1^6 h_6 e_0 g$ in degree $(106,15,60)$.  If there were a differential
\[
d_2 \left( \texttt{\{13-300\}} \right) = \texttt{\{15-313\}},
\]
then this differential would be detected by inclusion $S \rightarrow S/\tau$ of the bottom cell of $S/\tau$, since $\texttt{\{15-313\}}$ is not divisible by $\tau$.  However, according to \cite{IWX22}, no such differential occurs in the Adams spectral sequence for $S/\tau$.
\end{proof}

\begin{lem}\label{D1h1^2i1-d3}
$(107, 13, 60)$ 
$d_3 \left( \emph{\texttt{\{13-300\}}} \right) = 0$.
\end{lem}

\begin{proof}
By inspection, the only possible non-zero value for $d_3 \left( \texttt{\{13-300\}} \right)$ is $\texttt{\{16-344\}}$, also known as $h_2 g^2 C'$ in degree $(106,16,60)$.  We will show that $\tau \texttt{\{16-344\}}$ is non-zero on the $E_3$-page.  This will rule out the possible differential since $\texttt{\{13-300\}}$ is annihilated by $\tau$.

We know that $\tau \texttt{\{16-344\}}$ is non-zero in the $E_2$-page, so we just need to show that it cannot be hit by a $d_2$ differential.  The only possible source of such a differential is $\texttt{\{14-306\}}$ in degree $(107,14,59)$.  Note that $\texttt{\{14-306\}}$ equals $h_1 \texttt{\{13-293\}}$, but $\tau \texttt{\{16-344\}}$ is not a multiple of $h_1$.  Therefore, $d_2 \left( \texttt{\{14-306\}} \right)$ cannot equal
$\tau \texttt{\{16-344\}}$.
\end{proof}

There are no possible non-zero values for $d_4 \left( \texttt{\{13-300\}} \right)$ or $d_5 \left( \texttt{\{13-300\}} \right)$.

\begin{lem}
\label{D1h1^2i1-d6-partial}
$(107, 13, 60)$
$d_6 \left( \emph{\texttt{\{13-300\}}} \right)$ does not equal $\emph{\texttt{\{19-447\}}}$ nor $\emph{\texttt{\{19-447\}}} + \emph{\texttt{\{19-449\}}}$.
\end{lem}

\begin{proof}
Inclusion $S \rightarrow S/\tau$ of the bottom cell takes the element $\texttt{\{19-447\}}$ to $h_1^5 x_{101,14}$ on $E_2$-pages.  The latter element is still non-zero on the $E_6$-page of the Adams spectral sequence for $S/\tau$, and it is not hit by a $d_6$ differential.  Therefore, $\texttt{\{19-447\}}$ cannot be hit by a $d_6$ differential in the Adams spectral sequence for $S$.

Similarly, inclusion of the bottom cell takes $\texttt{\{19-447\}} + \texttt{\{19-449\}}$ to $h_1^5 x_{101,14} + \tau M h_1 g^3$.  The latter term is hit by a $d_2$ differential in the Adams spectral sequence for $S/\tau$, so inclusion of the bottom cell takes $\texttt{\{19-447\}} + \texttt{\{19-449\}}$ to $h_1^5 x_{101,14}$ on Adams $E_3$-pages.  The rest of the argument is identical to the previous paragraph.
\end{proof}

Our next task is to rule out the possibility that $d_6 \left( \texttt{\{13-300\}} \right)$ might equal $\texttt{\{19-449\}}$.  This case is somewhat harder than the cases considered in \cref{D1h1^2i1-d6-partial}, and it requires further study of related elements.

\begin{lem}
\label{17-398}
$(107, 17, 59)$
$d_2 \left( \emph{\texttt{\{17-398\}}} \right) = 0$.
\end{lem}

\begin{proof}
The only possible non-zero value is $\tau \texttt{\{19-449\}}$.  However, $\texttt{\{17-398\}}$ is a multiple of $h_3$, while $\tau \texttt{\{19-449\}}$ is not.
\end{proof}

\begin{lem}
\label{17-339}
$(100, 17, 56)$
$d_2(\emph{\texttt{\{17-339\}}}) = \emph{\texttt{\{19-393\}}}$.
\end{lem}

\begin{proof}
This follows by comparison along the inclusion $S \rightarrow S/\tau$ of the bottom cell of $S/\tau$.  The $d_2$ differential is present in the Adams spectral sequence for $S/\tau$ \cite{IWX22}.
\end{proof}

\begin{lem}
\label{14-268}
$(101, 14, 56)$
The element $\emph{\texttt{\{14-268\}}}$ survives to the $E_4$-page.
\end{lem}

\begin{proof}
There are no possible non-zero values for $d_2( \texttt{\{14-268\}} )$.  There are also no possible non-zero values for $d_3( \texttt{\{14-268\}} )$ because $\texttt{\{17-339\}}$ does not survive to the $E_3$-page by \cref{17-339}.
\end{proof}

\begin{lem}
\label{16-356}
$(107, 16, 60)$
The element $\emph{\texttt{\{16-356\}}}$ survives to the $E_4$-page.
\end{lem}

\begin{proof}
Note that $\texttt{\{16-356\}}$ equals $h_2^2 \texttt{\{14-268\}}$, so the result follows from \cref{14-268}.
\end{proof}

\begin{lem}
\label{13-288}
$(105, 13, 57)$
The element $\emph{\texttt{\{13-288\}}}$ survives to the $E_5$-page.
\end{lem}

\begin{proof}
The only possible non-zero value for $d_2( \texttt{\{13-288\}})$ is $\tau \texttt{\{15-293\}}$.  However, $h_2 \cdot \tau \texttt{\{15-293\}}$ is non-zero, while $h_2 \texttt{\{13-288\}}$ is zero.

There are no possible non-zero values for $d_3( \texttt{\{13-288\}})$ or $d_4( \texttt{\{13-288\}})$.
\end{proof}

\begin{lem}
\label{15-323}
$(107, 15, 59)$
The element $\emph{\texttt{\{15-323\}}}$ survives to the $E_5$-page.
\end{lem}

\begin{proof}
This follows from \cref{13-288} and the relation $\texttt{\{15-323\}} = h_1^2 \texttt{\{13-288\}}$.
\end{proof}

\begin{lem}
\label{14-306}
$(107, 14, 59)$
$d_2 \left( \emph{\texttt{\{14-306\}}} \right) = \emph{\texttt{\{16-352\}}}$.
\end{lem}

\begin{proof}
This follows by comparison along the inclusion $S \rightarrow S/\tau$ of the bottom cell of $S/\tau$.  The $d_2$ differential is present in the Adams spectral sequence for $S/\tau$ \cite{IWX22}.

This shows that $d_2 \left( \texttt{\{14-306\}} \right)$ equals either $\texttt{\{16-352\}}$ or $\texttt{\{16-352\}} + \tau \texttt{\{16-344\}}$.  The element $\texttt{\{14-306\}}$ is a multiple of $h_1$, but $\texttt{\{16-352\}} + \tau \texttt{\{16-344\}}$ is not a multiple of $h_1$.  This rules out the latter possibility.
\end{proof}

\begin{lem}
\label{20-485}
$(106, 20, 60)$
$d_2 \left( \emph{\texttt{\{20-485\}}} \right) = \emph{\texttt{\{22-548\}}}$.
\end{lem}

\begin{proof}
This follows by comparison along the inclusion $S \rightarrow S/\tau$ of the bottom cell.  The differential is present in the Adams spectral sequence for $S/\tau$ \cite{IWX22}.
\end{proof}

\begin{lem}
\label{17-391}
$(107, 17, 61)$
$d_2 \left( \emph{\texttt{\{17-391\}}} \right) = 0$, and $d_3 \left( \emph{\texttt{\{17-391\}}} \right) = \emph{\texttt{\{20-480\}}}$.
\end{lem}

\begin{proof}
Both differentials are detected by comparison along the inclusion $S \rightarrow S/\tau$ of the bottom cell.  The only possible non-zero value for $d_2 \left( \texttt{\{17-391\}} \right)$ is $\texttt{\{19-448\}}$.  If such a differential occurred, then it would be detected by a $d_2$ differential in the Adams spectral sequence for $S/\tau$.  However, there is no such differential for $S/\tau$.

The only possible non-zero value for $d_3 \left( \texttt{\{17-391\}} \right)$ is $\texttt{\{20-480\}}$.  In this case, there is a corresponding $d_3$ differential in the Adams spectral sequence for $S/\tau$.
\end{proof}

\begin{lem}
\label{19-449}
$(106,19,59)$ 
The element $\tau \emph{\texttt{\{19-449\}}}$ is non-zero in the Adams $E_6$-page.
\end{lem}

\begin{proof}
By inspection of the $E_2$-page, there are several possible differentials that might kill $\tau \texttt{\{19-449\}}$.  The table shows each possibility and gives the reason why that differential does not in fact occur.

\begin{longtable}{ll}
\caption{$\C$-motivic Adams differentials that could possibly equal $\tau \texttt{\{19-449\}}$
\label{table:diffcolor}
} \\
\toprule
differential & proof \\
\midrule \endfirsthead
\caption[]{$\C$-motivic Adams differentials that could possibly equal $\tau \texttt{\{19-449\}}$}\\
\toprule
differential & proof \\
\midrule \endhead
\bottomrule \endfoot
$d_2 \left( \texttt{\{17-398\}} \right)$ & \cref{17-398} \\
$d_2 \left( \tau^2 \texttt{\{17-391\}} \right)$ & \cref{17-391} \\
$d_3 \left( \tau \texttt{\{16-356\}} \right)$ & \cref{16-356} \\
$d_4 \left( \texttt{\{15-323\}} \right)$ & \cref{15-323} \\
$d_5 \left( \texttt{\{14-306\}} \right)$ & \cref{14-306} \\
\end{longtable}
\end{proof}

\begin{lem}
\label{D1h1^2i1-d6}
$(107, 13, 60)$
$d_6 \left( \emph{\texttt{\{13-300\}}} \right) = 0$.
\end{lem}

\begin{proof}
Because of \cref{D1h1^2i1-d6-partial}, it only remains to show that $d_6 \left( \texttt{\{13-300\}} \right)$ cannot equal $\texttt{\{19-449\}}$.  \cref{19-449} shows that $\tau \texttt{\{19-449\}}$ is non-zero in the $E_6$-page, while $\tau \texttt{\{13-300\}}$ is zero.
\end{proof}

\begin{lem}
\label{D1h1^2i1-d7}
$(107, 13, 60)$
$d_7 \left( \emph{\texttt{\{13-300\}}} \right) = 0$.
\end{lem}

\begin{proof}
In the $E_2$-page, there are two possible non-zero values $\texttt{\{20-485\}}$ and $\tau \texttt{\{20-480\}}$ for the differential.  However, \cref{20-485} and \cref{17-391} show that these elements do not occur in the $E_7$-page.
\end{proof}

\begin{prop}
\label{prop:13-300}
$(107, 13, 60)$
The element $\emph{\texttt{\{13-300\}}}$ is a permanent cycle.
\end{prop}

\begin{proof}
The previous lemmas have shown that $\texttt{\{13-300\}}$ survives to the $E_8$-page.  For higher differentials, all possible non-zero values disappear already in the $E_5$-page.
\end{proof}

\begin{prop}
\label{Delta1h1^2i1g^4k}
The element $\Delta_1 h_1^2 i_1 g^{4k}$ is a permanent cycle in the $\C$-motivic Adams spectral sequence for all $k \geq 0$.
\end{prop}

\begin{proof}
\cref{prop:no-differentials} shows that $\Delta_1 h_1^2 i_1 g^{4k}$ is non-zero if it does not support any differentials.

The proof proceeds by induction.  The base case $k=0$ says that $\Delta_1 h_1^2 i_1$ is a permanent cycle, which we know from \cref{prop:13-300}.
Let $\alpha_0$ be an element of $\pi_{107,60}$ that is detected by $\Delta_1 h_1^2 i_1$.  Note that $\eta^{15} \alpha_0$ is zero by part (2) of \cref{thm:C-homotopy}.

We inductively define $\alpha_k$ to be an element of the  Toda bracket $\langle \alpha_{k-1}, \eta^{15}, \eta_6 \rangle$.  Part (2) of \cref{thm:C-homotopy} ensures that $\alpha_{k-1} \cdot \eta^{12}$ is zero, so the bracket is well-defined.

The induction hypothesis is that $\alpha_{k-1}$ is detected by $\Delta_1 h_1^2 i_1 g^{4k-4}$.  The Moss convergence theorem \cite[Theorem 1.1]{BK21} implies that the Toda bracket $\langle \alpha_{k-1}, \eta^{15}, \eta_6 \rangle$ is detected by the Massey product $\langle \Delta_1 h_1^2 i_1 g^{4k-4}, h_1^{15}, h_1 h_6 \rangle$, which  equals $\Delta_1 h_1^2 i_1 g^{4k}$ as observed in \cref{tab:Massey}.  The Moss convergence theorem requires that there be no ``crossing differentials''.  In our case, such crossing differentials would have sources with negative Chow degree, but there are no such elements by \cref{prop:Chow-degree-nonneg}.  Therefore, $\alpha_k$ is detected by $\Delta_1 h_1^2 i_1 g^{4k}$.
\end{proof}

\section{The Adams--Novikov \texorpdfstring{$E_2$}{E2}-page}
\label{sctn:Adams-Novikov}

The purpose of this section is to explain how \cref{cor:Adams-Novikov-E2}, \cref{thm:Adams-Novikov-E2}, and \cref{thm:Adams-Novikov-E2-beta} follow from our main  results.  These corollaries translate our results about $\C$-motivic homotopy into results about the cohomology of the moduli stack of 1-dimensional formal group laws, also known as the classical Adams--Novikov $E_2$-page.

We recall some standard notation from the Adams--Novikov $E_2$-page:
\begin{itemize}
\item
We grade elements in the Adams--Novikov $E_2$-page in the form
$(s,f)$, where $s$ is the stem and $f$ is the Adams--Novikov filtration.
\item 
$\alpha_1$ is the element in degree $(1,1)$ that detects $\eta$.  It corresponds to $h_1$ in the Adams $E_2$-page.
\item 
$\beta_{2^j/2^j -1}$ is the ``beta family'' element in degree $(2^{j+2}, 2)$ that detects Mahowald's $\eta_{j+2}$.
\end{itemize}

We use the terminology and notation of \cref{subsctn:coexponent} to describe $\alpha_1$-submodules of the Adams--Novikov $E_2$-page in terms of cogenerators and $\alpha_1$-coexponents.  Beware that the Adams--Novikov $E_2$-page is not an $\F_2$-algebra.  We are primarily interested in elements that happen to be $2$-torsion in the Adams--Novikov $E_2$-page.  However, the generators of our $\alpha_1$-submodules may not be 2-torsion, although we do not know of any specific examples. The example involving $\eta_6$ discussed in \cref{subsctn:w1-periodic-homotopy} does not apply, since $\beta_{16/15}$ is $2$-torsion in the Adams--Novikov $E_2$-page, but it supports a hidden $2$ extension in the Adams--Novikov $E_\infty$-page.

\begin{cor}
\label{thm:Adams-Novikov-E2}
For all $n \geq 1$, there are non-zero elements $a_n$ in degree $(20n-2, 4n)$ in the Adams--Novikov $E_2$-page such that:
\begin{enumerate}
\item
$a_n$ is annihilated by $\alpha_1$.
\item
the elements form a $w_1$-periodic family in the sense that
$a_{n+1}$ is contained in the Massey product $\langle a_n, \alpha_1, a_1 \rangle$ for all $n$.
\item 
the element $a_{2^j}$ equals $\alpha_1^{2^{j+2} - 2} \cdot \beta_{2^{j+2}/2^{j+2} -1}$.
\item 
the $\alpha_1$-coexponent of $a_n$ is given  in \cref{table:AN-submodules} for some values of $n$.
\end{enumerate}
\end{cor}

\begin{proof}
According to \cref{prop:no-tau-extns-in}, the elements of the form $\frac{\overline{w_1^{4k}}}{h_1^n}$ in the Adams $E_\infty$-page detect elements in $\C$-motivic homotopy that are not divisible by $\tau$.  Therefore, they map under inclusion of the bottom cell to non-zero elements in the homotopy of the cofiber of $\tau$.  The homotopy of the cofiber of $\tau$ is isomorphic to the Adams--Novikov $E_2$-page \cite[Theorem 6.12]{Isa19}.  Parts (1), (2), and (4) are translations of \cref{thm:C-homotopy} along this isomorphism.  Part (3) follows from the fact that both the Adams $E_\infty$-page element $h_1 h_{j+4}$ and Adams--Novikov $E_\infty$-page element $\beta_{2^{j+2}/2^{j+2} -1}$ detect Mahowald's $\eta_{j+4}$.
\end{proof}

\begin{rem}
In part (4) of \cref{thm:Adams-Novikov-E2}, the $\alpha_1$-coexponent of $a_n$ is unknown for many values of $n$.  The smallest such value is $n = 12$.  See \cref{rem:unknown-coexponents} for further discussion of these uncertainties.
\end{rem}

\begin{cor}
\label{thm:Adams-Novikov-E2-beta}
For all $n \geq 1$, there are non-zero elements $b_n$ in degree $(20n+9, 4n+3)$ in the Adams--Novikov $E_2$-page such that:
\begin{enumerate}
\item
$b_n$ is annihilated by $\alpha_1$.
\item
the elements form a $w_1$-periodic family in the sense that $b_{n+1}$ is contained in the Massey product $\langle b_n, \alpha_1, a_1 \rangle$ for all $n$.
\item 
the $\alpha_1$-coexponent of $b_n$ is 2.
\end{enumerate}
\end{cor}

\begin{proof}
Similarly to the proof of \cref{thm:Adams-Novikov-E2}, these results are translations of \cref{thm:C-homotopy-beta} along the isomorphism between the Adams--Novikov $E_2$-page and the homotopy of the cofiber of $\tau$.
\end{proof}

\cref{table:AN-submodules} describes some $\alpha_1$-submodules of the Adams--Novikov $E_2$-page in terms of their cogenerators and $\alpha_1$-coexponents.  The generators do not have standard names.  Translations of the degrees are displayed because they reveal the numerical patterns more evidently.

\begin{longtable}{lll}
\caption{Some $\alpha_1$-submodules of the classical Adams--Novikov $E_2$-page
\label{table:AN-submodules}}
\\
\toprule
cogenerator & $\mathrm{degree} + (2,0)$ & $\alpha_1$-coexponent \\
\midrule \endfirsthead
\caption[]{Some $\alpha_1$-submodules of the classical Adams--Novikov $E_2$-page}\\
\toprule
cogenerator & $\mathrm{degree} + (2,1)$ & $\alpha_!$-coexponent \\
\midrule \endhead
\bottomrule \endfoot
$a_{2^j}$ & $2^j \cdot (20, 4)$ & $2^{j+2} - 1$ \\
$b_k$ & $k \cdot (20,4) + (11, 3)$ & $2$ \\
$a_{2k+3}$ & $(2k+3) \cdot (20, 4)$ & $6$ \\
$a_{4k+6}$ & $(4k+6) \cdot (20, 4)$ & $12$  \\
\end{longtable}

\begin{rem}
The $\C$-motivic homotopy elements $\gamma_n$ and $\delta_n$ of \cref{thm:C-homotopy} and \cref{thm:C-homotopy-beta} are detected respectively by the elements $a_n$ and $b_n$ from \cref{thm:Adams-Novikov-E2} and \cref{thm:Adams-Novikov-E2-beta} in the $\C$-motivic Adams--Novikov $E_\infty$-page.
\end{rem}

\section{The cohomology of \texorpdfstring{$\C$}{C}-motivic \texorpdfstring{$\cA(2)$}{A(2)} in Chow degree one}
\label{sctn:A(2)}

As an illustration of our methods with the Burklund--Xu spectral sequence from \cref{sctn:BXss}, we carry out the complete computation of the cohomology of $\C$-motivic $\cA(2)$ in Chow degree one.  This section can serve as a training module for readers who seek hands-on experience with the Burklund--Xu spectral sequence.  We highly recommend that such readers reconstruct the material in this guide for themselves.

The cohomology of $\C$-motivic $\cA(2)$ can be computed manually \cite{May-A(2)} \cite{Isa09}, but it is a lengthy computation whose structure cannot be fully appreciated at one glance.  On the other hand, if we restrict to Chow degrees zero and one, the computation becomes easily manageable.

\begin{prop}
\label{A(2)-Cd0}
When restricted to Chow degree zero, the cohomology of $\C$-motivic $\cA(2)$ is
\[
\frac{\mathbb{F}_2 [h_1, h_2, u, g]}
{h_1 h_2, h_2^3, h_2 u, u^2 + h_1^2 g}.
\]
\end{prop}

\begin{proof}
The dual $\cA(2)_*$ of $\cA(2)$ is
\[
\frac{\mathbb{F}_2[\tau] [\tau_0, \tau_1, \tau_2, \xi_1, \xi_2]}
{\tau_0^2 + \tau \xi_1, \xi_1^4, \tau_1^2 + \tau \xi_2, \xi_2^2,
\tau_2^2}
\]
Define $\cB(2)_*$ to be the subalgebra
\[
\frac{\mathbb{F}_2[\xi_1, \xi_2]}
{\xi_1^4, \xi_2^2}
\]
of $\cA(2)_*$, and note that the cobar complex of $\cB(2)$ is equal to the part of the cobar complex of $\cA(2)$ in Chow degree zero.  Therefore, $H^{***}\cB(2)$ is isomorphic to the part of $H^{***} \cA(2)$ in Chow degree zero.

It remains to compute $H^{***}\cB(2)$.  Observe that $\cB(2)_*$ is isomorphic to classical $\cA^{\cl}(1)_*$, with shifted degrees.  Consequently, we thoroughly know the cohomology of $\cB(2)$ by comparison to the well-known cohomology
\[
\frac{\mathbb{F}_2 [h_0, h_1, a, P]}{h_0 h_1, h_1^3, h_1 a, a^2 + P h_0^2}
\]
of classical $\cA^{\cl}(1)$. 
\end{proof}

Now we proceed to Chow degree one.

\begin{figure}[H]
\begin{center}
\includegraphics[trim={0cm, 0cm, 0cm, 0cm},clip,page=1,scale=0.55]{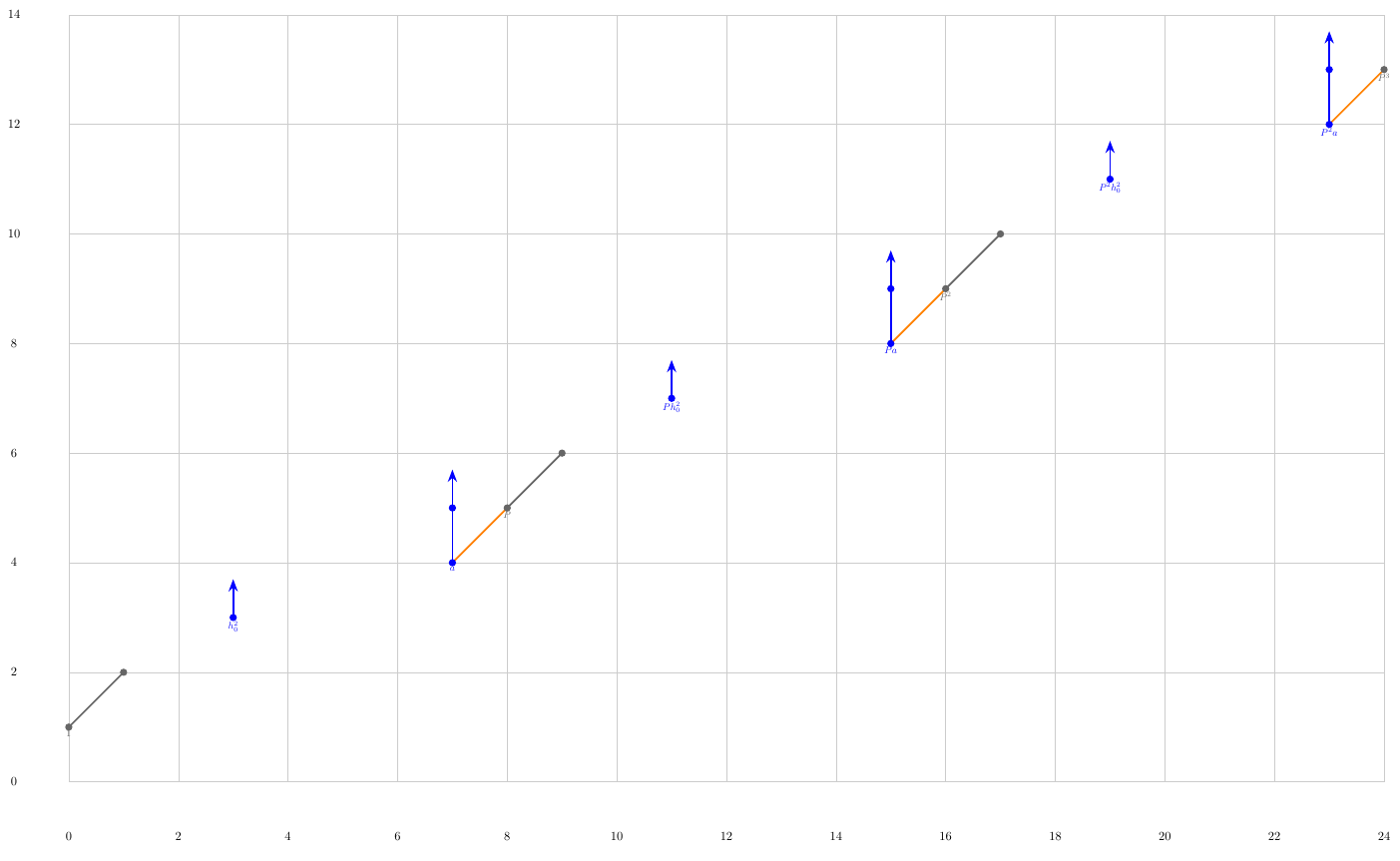}
\caption{The $E_{\infty}$-page of the Burklund--Xu spectral sequence for the cohomology of $\C$-motivic $\cA(2)$ in Chow degree one}
\label{fig:A(2)-Einfty}
\end{center}
\end{figure}

\begin{prop}
\label{A(2)-Cd1}
The $E_\infty$-page of the Burklund--Xu spectral sequence for the cohomology of $\C$-motivic $A(2)$ in Chow degree one is displayed in \cref{fig:A(2)-Einfty}.
\end{prop}

\begin{figure}[H]
\begin{center}
\includegraphics[trim={0cm, 0cm, 0cm, 0cm},clip,page=1,scale=0.55]{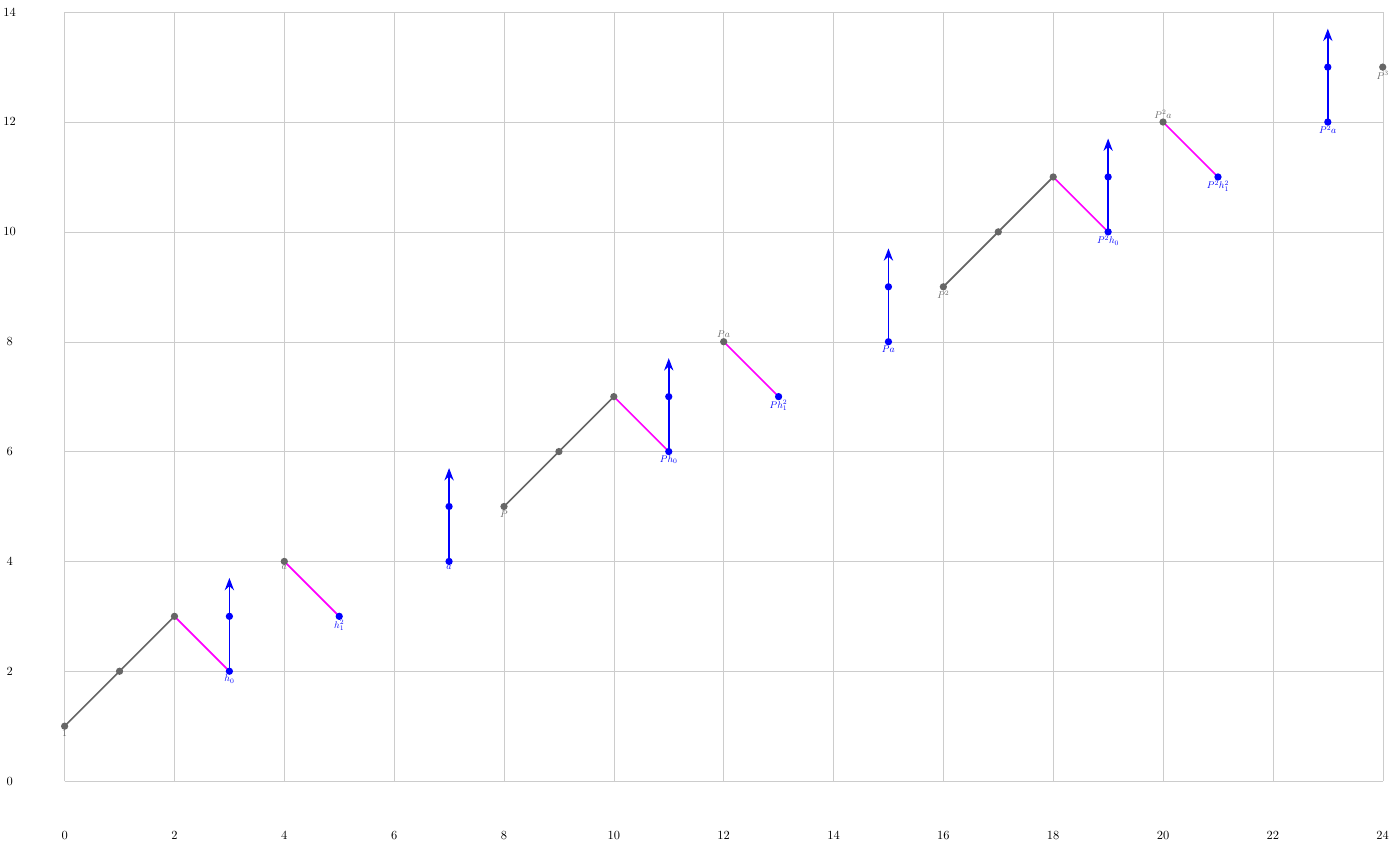}
\caption{The $E_2$-page of the Burklund--Xu spectral sequence for the cohomology of $\C$-motivic $\cA(2)$ in Chow degree one}
\label{fig:A(2)-E2}
\end{center}
\end{figure}

\begin{proof}
Because $\cA(2)_*$ contains $\tau_0$, $\tau_1$, and $\tau_2$, but no $\tau_i$ for $i \geq 3$, the $E_1$-page of the $\C$-motivic Burklund--Xu spectral sequence in Chow degree one takes the form
\[
q_0 \cdot H^{**}\cA^{\cl}(1) \oplus 
q_1 \cdot H^{**}\cA^{\cl}(1) \oplus 
q_2 \cdot H^{**}\cA^{\cl}(1)
\]
with differentials $d_1(q_1) = q_0 \cdot h_0$ and $d_1(q_2) = q_1 \cdot h_1$.  Beware of the notation that we are using for the coefficients on $q_n$, as discussed carefully in \cref{sctn:BX-Cd1}.

The $E_2$-page now takes the form
\[
q_0 \cdot \coker(h_0) \bigoplus 
q_1 \cdot \frac{\ker(h_0)}{\im(h_1)} \bigoplus 
q_2 \cdot \ker(h_1),
\]
where the kernels and images are taken within $H^{**}\cA^{\cl}(1)$.  More explicitly, the $E_2$-page is
\[
q_0 \cdot \frac{\mathbb{F}_2 [h_1, a, P]}{h_1^3, h_1 a, a^2}  \bigoplus
q_1 \cdot 0 \bigoplus
q_2 \cdot 
\left( h_0 \cdot \mathbb{F}_2 [h_0, P] \oplus
h_1^2 \cdot \mathbb{F}_2 [P] \oplus
a \cdot \mathbb{F}_2[h_0, P] \right).
\]
Note the curious fact that $\frac{\ker(h_0)}{\im(h_1)}$ equals $0$ in $H^{**}\cA^{\cl}(1)$, so the $q_1$-summand contributes nothing to the $E_2$-page.  \cref{fig:A(2)-E2} illustrates the $E_2$-page that we have just described symbolically.

Next, we must analyze $d_2$ differentials, which can only take the form $d_2(q_2 \cdot x) = q_0 \cdot y$.  For degree reasons, there are a few but not many possible differentials.  Namely, we have
\begin{gather*}
d_2(q_2 \cdot P^k h_0) = q_0 \cdot P^k h_1^2 \\
d_2(q_2 \cdot P^k h_1^2 ) = q_0 \cdot P^k a.
\end{gather*}
As in \cref{rem:higher-diff}, these formulas correspond to Massey product computations
\begin{gather*}
P^k h_1^2 = \langle P^k h_0, h_1, h_0 \rangle \\
P^k a = \langle P^k h_1^2, h_1, h_0 \rangle.
\end{gather*}

There are no further differentials, and we directly obtain the $E_3$-page, which equals the $E_\infty$-page.
\end{proof}

We provide a key for reading \cref{fig:A(2)-E2} and \cref{fig:A(2)-Einfty}. Dots represent copies of $\mathbb{F}_2$.  Multiples of $q_0$ appear in black, and multiples of $q_2$ appear in blue. (There are no multiples of $q_1$, as explained in the proof of \cref{A(2)-Cd1}.)  The labels indicate the coefficients on $q_n$ for all elements.

As usual, vertical lines indicate multiplications by $h_0$, and lines of slope $1$ indicate multiplications by $h_1$.  Vertical arrows indicate infinitely many $h_0$ multiplications.  (Because of our notational conventions, these lines represent multiplications by $h_1$ and $h_2$ respectively in the cohomology of $\C$-motivic $\cA(2)$.)

In \cref{fig:A(2)-E2}, the $d_2$ differentials appear as magenta lines of slope $-1$.  In \cref{fig:A(2)-Einfty}, hidden $h_1$ extensions appear as orange lines.

\begin{lem}
\label{A(2)-hidden}
In the $\C$-motivic Burklund--Xu spectral sequence for the cohomology of $A(2)$ in Chow degree one, there are no hidden $h_0$ extensions.  The only hidden $h_1$ extensions are from $q_2 \cdot P^k a$ to $q_0 \cdot P^{k+1}$.
\end{lem}

\begin{proof}
For degree reasons, there are no possible hidden $h_0$ extensions, and the only possible hidden $h_1$ extensions are from $q_2 \cdot P^k a$ to $q_0 \cdot P^{k+1}$.

Recall the Massey product $\langle a, h_1, h_1^2 \rangle = P h_0$ in $H^{**} \cA^{\cl}(1)$.  Then there is a similar Massey product $\langle q_2 \cdot a, h_1, h_1^2 \rangle = q_2 \cdot P h_0$ in the Burklund--Xu $E_1$-page.  If this Massey product were well-defined in $H^{***} \cA(2)$, then the May convergence theorem \cite{May69} would imply that $q_2 \cdot P h_0$ is a permanent cycle in the Burklund--Xu spectral sequence.  But we already saw in the proof of \cref{A(2)-Cd1} that $q_2 \cdot P h_0$ is not a permanent cycle.  Therefore, the Massey product is not well-defined in $H^{***} \cA(2)$.  The only possible problem is that there is a hidden $h_1$ extension on $q_2 \cdot a$, and $q_0 \cdot P$ is the only possible target for this hidden extension.

The rest of the hidden extensions follow from the first one by $P^k$ multiplication.
\end{proof}

\begin{rem}
Converting from our Burklund--Xu spectral sequence notation to standard $\C$-motivic notation for $H^{***} \cA(2)$, the $E_\infty$-elements $q_2 \cdot a$ and $q_0 \cdot P$ detect $e_0$ and $h_0 g$ respectively.  The hidden $h_1$ extension from $q_2 \cdot a$ to $q_0 \cdot P$ corresponds to the relation $h_2 \cdot e_0 = h_0 g$. This relation is elementary from the perspective of the $\C$-motivic May spectral sequence \cite{Isa09}.
\end{rem}

\cref{table:A(2)-Cd1} explicitly list the elements in $H^{***} \cA(2)$ in Chow degree one.  All of these elements are $g$-free, in the sense that $x \cdot g^k$ is non-zero and has Chow degree $1$ for every $x$ in the table and every $k \geq 0$.  The third column gives representatives in the Burklund--Xu $E_\infty$-page.  The fourth column indicates that some of the elements are $h_1$-periodic, in the sense that they support infinitely many multiplications by $h_1$.

\begin{longtable}{llll}
\caption{The cohomology of $\C$-motivic $\cA(2)$ in Chow degree one, up to $g$-multiples
\label{table:A(2)-Cd1}
} \\
\toprule
$(s,f,s-w)$ & element & detected by & $h_1$-periodic \\
\midrule \endfirsthead
\caption[]{The cohomology of $\C$-motivic $\cA(2)$ in Chow degree one, up to $g$-multiples}\\
\toprule
$(s,f,s-w)$ & element & detected by & $h_1$-periodic \\
\midrule \endhead
\bottomrule \endfoot
$(0,1,0)$ & $h_0$ & $q_0 \cdot 1$ \\
$(3,1,1)$ & $h_0 h_2$ & $q_0 \cdot h_1$ \\
$(8,3,3)$ & $c_0$ & $q_2 \cdot h_0^2$ & $h_1$-periodic \\
$(17,4,7)$ & $e_0$ & $q_2 \cdot a$ & $h_1$-periodic
\end{longtable}

\bibliographystyle{amsalpha}
\bibliography{bib}
\end{document}